\documentclass[leqno,12pt]{article}
\usepackage{a4}
\usepackage{amsthm,amssymb}
\usepackage{amsbsy,amsfonts,amsmath}

\usepackage{latexsym}

\numberwithin{equation}{section}

\newtheorem{lmm}{Lemma}[section]

\newtheorem{thm}{Theorem}[section]

\theoremstyle{definition}
\newtheorem{dfn}{Definition}[section]
\theoremstyle{remark}
\newtheorem{rem}{Remark}[section]%%%
\def\al{\alpha}

\def\kbka{{\frac {k(\beta)} {k(\al)}}}

\def\a{{\mathcal A}}
\def\b{{\mathcal B}}
\def\bp{{{\mathcal B}_+}}
\def\btp{{{\mathcal B}^{2,\prime}}}

\def\r{\mathcal R}

\def\v{\mathcal V}
\def\vo{{\mathcal V}^0}
\def\vt{{\mathcal V}^\times}

\def\vc{\v^{\mathbb C}}

\def\bbbc{{\mathbb C}}
\def\bbbct{{\mathbb C}^\times}

\def\bbbf{{\mathbb F}}
\def\bbbn{{\mathbb N}}

\def\bbbr{{\mathbb R}}
\def\bbbrt{{\mathbb R}^\times}
\def\bbbx{{\mathbb X}}
\def\bbbxp{{\mathbb X}_+}
\def\bbbxm{{\mathbb X}_-}
\def\bbbz{{\mathbb Z}}
\def\bbbzp{{\mathbb Z}_+}
\def\bbbzm{{\mathbb Z}_-}
\def\bbbzt{{\mathbb Z}^\times}
\def\fa{\mathfrak{a}}
\def\fb{\mathfrak{b}}
\def\fg{\mathfrak{g}}
\def\fh{\mathfrak{h}}
\def\fn{\mathfrak{n}}

\def\ad{{\rm{ad}}}

\def\Rsh{R_{{\rm{sh}}}}
\def\Rlg{R_{{\rm{lg}}}}
\def\Rex{R_{{\rm{ex}}}}

\def\Lsh{L_{{\rm{sh}}}}
\def\Llg{L_{{\rm{lg}}}}
\def\Lex{L_{{\rm{ex}}}}

\def\PIsh{\Pi_{{\rm{sh}}}}
\def\PIlg{\Pi_{{\rm{lg}}}}
\def\PIex{\Pi_{{\rm{ex}}}}

\providecommand{\bysame}{\leavevmode\hbox to3em{\hrulefill}\thinspace}
\providecommand{\MR}{\relax\ifhmode\unskip\space\fi MR }
\providecommand{\href}[2]{#2}

\title{Exposition on affine and elliptic root systems\\ and elliptic Lie algebras}

\author{Saeid Azam,\quad Hiroyuki Yamane, \quad Malihe Yousofzadeh}

\date{}

\begin{document}
\maketitle

%\tableofcontents      %optional

\begin{abstract}
This is an exposition in order to give an explicit way to understand
(1) a non-topological proof for an existence of a base of an affine root system, 
(2) a Serre-type definition of an elliptic Lie algebra with rank $\geq 2$, 
and (3) the isotropic root multiplicities  
 obtained  
 from a viewpoint of the Saito-marking lines. 
\end{abstract}

\section{Introduction}
In 1985, K.~Saito~\cite{saito1} introduced the notion of an
{\it{$n$-extended affine root system}}. If $n=0$
(respectively, $n=1$), it is an irreducible finite root
system
(respectively, an affine root system).
In~\cite{saito1}, he also intensively
studied $2$-extended affine root systems, which are now
called {\it{elliptic root systems}} (see
\cite{saitotakebayashi}).
Since then, various attempts have been made to
construct Lie algebras whose non-isotropic roots
form those root systems. Among them are
{\it{toroidal Lie algebras}} \cite{MoodyRaoYokonuma},
{\it{extended affine Lie algebras}}
\cite{Allisonetal},
and {\it{toral type extended affine Lie algebras}} \cite{MoodyRaoYokonuma}, 
\cite{Yousofzadeh}.
See \cite[Introduction]{saitoyoshii} for the history. 

\par In 2000, K.~Saito and D.~Yoshii~\cite{saitoyoshii}
constructed certain Lie 
algebras by
using the Borcherds
lattice vertex algebras, called them
{\it{simply-laced elliptic Lie algebras}} and showed that
they are
isomorphic
to $ADE$-type ($2$-variable) toroidal Lie algebras of
rank $\geq 2$.
They also gave two other definitions for their Lie
algebras. One uses
generators and relations. The other uses (affine-type)
Heisenberg Lie algebras; this was generalized by
D.~Yoshii~\cite{yoshii}
in order to define Lie algebras associated with the
reduced elliptic root systems,
and he called them {\it{elliptic Lie algebras}}. In 2004,
the second author~\cite{yamane2} gave defining relations of the
elliptic Lie algebras of rank $\geq 2$.

The aim of this paper is to give an exposition in order to give an explicit way to understand
(1) a non-topological proof for the existence of a base of an affine root system
(Theorem~\ref{theorem:theoremoneone}, originally given in \cite{MacDonald}), 
(2) a Serre-type definition of an elliptic Lie algebra $\fg$ with rank $\geq 2$
(Definition~\ref{definition:defofellsuper}, originally given in \cite{yamane2}) and 
the fact that the non-isotropic roots form the corresponding elliptic root system
and their multiplicities are one (Theorem~\ref{theorem:mmain}, originally given in \cite{yamane2}),  
and (3) a list of the multiplicities of the isotropic roots of $\fg$,  
proved from a viewpoint of the Saito-marking lines 
(Theorem~\ref{theorem:dimnull}, new result). 
 
As for (2), 
we point out that
our defining relations are closely related to
defining relations, called
{\it{Drinfeld realization}}, of the quantum affine
algebras
due to V.G.~Drinfeld~\cite[Theorems 3 and 4]{Drinfeld}.
Recently the same authors have written a paper \cite{AYY}, motivated by \cite{YouRIMS}, 
giving a finite 
number of defining relations of the universal coverings of some Lie tori.  

We hope that the material presented here regarding affine root systems, in particular     
the existence of a base, would give another point of view to readers interested in the subject, specially to those reading the book~\cite{MacDonald2} by I.G.~MacDonald. 
(Incidentally, in order to read \cite{MacDonald2}, 
we also hope that the paper \cite{HecYam08b} would also be 
helpful in being familiar with Coxeter groups, especially the Matsumoto theorem.)

\section{Preliminary}

In this section, we mention elemental properties of (Saito's) extended affine root systems,
especially  
(\ref{eqn:eqnbsofexaff}).

\subsection{Basic notation and terminology}\label{subsection:bcnt}\label{subsection:basicntn}
As usual, we let $\bbbz$ denote the ring of integers,
$\bbbn$ the set of positive integers, $\bbbr$
the field of real numbers, and
$\bbbc$ the field of complex numbers.
For a set $S$, let $|S|$ denote 
the cardinal number of $S$. If $S$ is a subset of $\bbbr$,
we let $S^\times:=\{s\in S|s\ne 0\}$, $S_+:=\{s\in S|s\geq 0\}$, and
$S_-:=\{s\in S|s\leq 0\}$. 
\par
For a unital subring $X$ of $\bbbc$, an $X$-module $M$, a subset $Y$ of $X$,
subsets $S$ and $S^\prime$ of $M$, $x\in X$
and $m\in M$, we let $S+S^\prime:=\{m+m^\prime\in M|m\in S,m^\prime \in S^\prime\}$,
$m+S:=\{m\}+S$,
$YS:=\{y_1s_1+\cdots+y_rs_r|r\in\bbbn,\,y_i\in Y,\,s_i\in S\,(1\leq i\leq r)\}$, 
$Ym:=Y\{m\}$, 
$xS:=\{x\}S$ and $-S:=(-1)S$; we understand $S+\emptyset=\emptyset$, $\emptyset S=\emptyset$
and $Y\emptyset =\emptyset$. \par
Throughout this paper,
for any $\bbbf$-linear space $\v$ with
a symmetric bilinear form $(\,,\,):\v\times\v
\rightarrow\bbbf$,
where $\bbbf$ is $\bbbr$ or $\bbbc$,
we set $\vo:=\{v\in\v|(v,v)=0\}$
and $\vt:=\v\setminus\vo$;
for each
$v\in\vt$, we  set
$v^\vee:={\frac {2v} {(v,v)}}$
and define $s_v\in{\rm GL}(\v)$ by
$s_v(z)=z-(v^\vee,z)v$
($z\in \v$); for any non-empty subset $S$
of $\vt$, we denote by $W_S$
the subgroup of ${\rm GL}(\v)$ generated
by $\{s_v|v\in S\}$, i.e.,
\begin{equation}\label{eqn:eqndefofws}
W_S:=\langle s_v|v\in S \rangle,
\end{equation} and moreover, 
let $W_S\cdot S^\prime:=\{w(z^\prime)\in\v |w\in W_S,z^\prime\in S^\prime\}$,
$W_S\cdot z:=W_S\cdot \{z\}$
for a subset $S^\prime$ of $\v$ and $z\in \v$,
and say that a subset $S$ of $\vt$ is {\it{connected}} if there exists no non-empty proper
subset $S^\prime$ of $S$ with $(S^\prime,S\setminus S^\prime)=\{0\}$.
For a subset $\v^\prime$ of $\v$, let $(\v^\prime)^0:=\v^\prime\cap\vo$,
and  $(\v^\prime)^\times:=\v^\prime\cap\vt$. We call an element of $\vo$ {\it{isotropic}}.

In this paper, we always
let
\begin{equation}\label{eqn:eqnpipr}
\pi:\v \rightarrow \v /\vo 
\end{equation} denote the canonical map. 

\subsection{Extended affine root systems}
\label{subsection:subsectiononeone}

\begin{dfn} \label{definition:definitionexaff}
{\rm{Let $l\in\bbbn$ and $n\in\bbbzp$.
Let $\v$ be an $(l+n)$-dimensional $\bbbr$-linear
space. Recall $\vo$ and $\vt$ from
Subsection~\ref{subsection:bcnt}.
Assume that there exists a positive
semi-definite symmetric bilinear form
$(\,,\,):\v\times \v\rightarrow\bbbr$
such that $\dim_\bbbr \vo =n$.
Let $R$ be a subset of $\v$. Then $R$
(or more precisely, $(R,\v)$) is an
{\it {{\rm{(}}$n$-{\rm{)}}extended affine root system}}
if $R$ satisfies
the following axioms:
\newline\par
({\rm{AX}}1) $R\subset \vt $, $\v =\bbbr R$. \par
({\rm{AX}}2) $\bbbz R$ is free as a $\bbbz$-module, and 
${\rm
rank}_\bbbz \bbbz R=n+l(=\dim_\bbbr \v )$. \par
({\rm{AX}}3) $(\al^\vee,\beta)\in\bbbz $ for $\al$,
$\beta\in R$. \par
({\rm{AX}}4) $s_\al(R)=R$ for all $\al\in R$. \par
({\rm{AX}}5) $R$ is connected.
%If $R=R_1\cup R_2$
%and $(\al,\beta)=0$ for $\al\in R_1$ and $\beta\in R_2$,
%then $R_1=\emptyset$ or $R_2=\emptyset$.
\newline\newline
(see \cite[(1.2) Definition~1 and (1.3) Note~2
iii)]{saito1} and see \cite{Azam} for an equivalence to 
\cite[Definition~2.1]{Allisonetal}.) 
%For $\mu\in R$, let
%\begin{equation}\label{eqn:eqnparitydash}
%p(\mu)=\left\{\begin{array}{ll}
%0 & \quad\mbox{if $2\mu\notin R$,} \\
%1 & \quad\mbox{if $2\mu\in R$}
%\end{array}\right.
%\end{equation} (see (\ref{eqn:eqnoneone}) below).
Let $W=W_R$ (see (\ref{eqn:eqndefofws})).}}
\end{dfn}
Let $R$ be as in Definition~\ref{definition:definitionexaff}.
It is well-known that for all
$\al\in R$,
\begin{equation}\label{eqn:eqnoneone}
\left\{\begin{array}{l}
R\cap\bbbr \al=\{\al, -\al\},\,
\{\al,2\al, -\al,-2\al\}\,\,{\rm {or}}\,\,
\{\al,{\frac {1} {2}}\al, -\al,-{\frac {1} {2}}\al\}, \\
\mbox{(so $-R=R$).}
\end{array}\right.
\end{equation} 
We call $R$ {\it {reduced}} (resp.
{\it {non-reduced}}) if $R\cap 2R=\emptyset$
(resp.
$R\cap 2R\ne\emptyset$).
\par
We say that two extended affine root systems $(R,\v )$
and $(R^\prime, \v^\prime)$ are {\it isomorphic}
if there exist an $\bbbr$-linear bijective map
$f:\v\rightarrow \v^\prime$ and $c\in\bbbr$
with $c>0$
such that
$f(R)=R^\prime$ and $(f(v),f(w))=c(v,w)$
for $v$, $w\in \v$. 
\begin{equation}\label{eqn:rtiso}
\mbox{We call this $f$ a {\it{root system isomomorphism}}.}
\end{equation} 

Let $R$, $l$ and $n$ be as above.
\par
By \cite[Theorem~5 of Chapter~XV]{lang}, 
since $\bbbz R/(\bbbz R)^0$ is torsion free, ({\rm{AX}}1-5) imply that there exists an $\bbbr$-basis
$\{x_1,\ldots,x_{l+n}\}$ of $\v$ such that
$\{x_{l+1},\ldots,x_{l+n}\}$ is an $\bbbr$-basis of $\vo$, $\{x_1,\ldots,x_{l+n}\}$
is a $\bbbz$-basis of the (torsion) free $\bbbz$-module $\bbbz R$
and $\{x_{l+1},\ldots,x_{l+n}\}$ is a $\bbbz$-basis of the (torsion) free $\bbbz$-module
$(\bbbz R)^0$ (see Subsection~\ref{subsection:basicntn} for notation), that is, 
\begin{equation}\label{eqn:eqnbsofexaff}
\left\{\begin{array}{l}
\v=\bbbr R=\bigoplus_{i=1}^{l+n}\bbbr x_i,\,\,\vo=\bigoplus_{j=l+1}^{l+n}\bbbr x_j, \\
\bbbz R=\bigoplus_{i=1}^{l+n}\bbbz x_i,\,\,(\bbbz R)^0=\bigoplus_{j=l+1}^{l+n}\bbbz x_j,
 \\
\dim_\bbbr \v={\rm{rank}}_\bbbz \bbbz R=n+l,
\,\,\dim_\bbbr \vo={\rm{rank}}_\bbbz (\bbbz R)^0=n.
\end{array}\right.
\end{equation}

Let $\{a_1,\ldots,a_n\}$ be a $\bbbz$-basis of $(\bbbz R)^0$.
Then there exist  $x_1,\ldots,x_l\in \bbbz R$
such that $\{x_1,\ldots,x_l,a_1,\ldots,a_n\}$ is 
a $\bbbz$-basis of $\bbbz R$ as well as  an $\bbbr$-basis of $\v=\bbbr R$
(see above). 
Let $1\leq m\leq n$. Let $\pi^\prime : \v \rightarrow \v /(\bbbr a_m\oplus\cdots \oplus\bbbr a_n)$ 
be the canonical map. Note that $\{\pi^\prime(x_1),\ldots,\pi^\prime(x_l),
\pi^\prime(a_1),\ldots,\pi^\prime(a_{m-1})\}$ is 
an $\bbbx$-basis of $\bbbx \pi^\prime(R)$ for $\bbbx\in\{\bbbz,\bbbr\}$.
In particular,
we see that
%\begin{equation}\label{eqn:bsofexaffd}
%\bbbz\pi^\prime(R)=\bigoplus_{i=1}^{l+m-1}\bbbz\pi^\prime(y_i)\Longrightarrow
%\bbbz R=(\bigoplus_{i=1}^{l+m-1}\bbbz y_i)\oplus (\bigoplus_{j=m}^n\bbbz a_j).
%\end{equation}
\begin{equation}\label{eqn:bsofexaffd}
\begin{array}{l}
\mbox{if $y_1,\ldots,y_{l+m-1}$ are elements of $\bbbz R$ such that} \\ 
\mbox{$\{\pi^\prime(y_1),\ldots,\pi^\prime(y_{l+m-1})\}$ is a $\bbbz$-base of the 
free $\bbbz$-module $\bbbz\pi^\prime(R)$,} \\
\mbox{then $\{y_1,\ldots,y_{l+m-1},a_m,\ldots,a_n\}$ is 
an $\bbbx$-basis of $\bbbx R$ for $\bbbx\in\{\bbbz,\bbbr\}$.}
\end{array} 
\end{equation}

\begin{equation}\label{eqn:ranknullitydef}
\mbox{We call $l$ the {\it {rank}}
of $R$.
We call $n$ the {\it {nullity}} of $R$.}
\end{equation}
If $n=0$, then $R$ is
an {\it {irreducible finite root system}}
(see \cite[(1.3) Example~1~i)]{saito1}).
If $n=1$, then $R$ is
an {\it {affine root system}} (see \cite[(1.3)
Example~1~ii)]{saito1}), see also Remark~\ref{remark:MacSaito}
below.
If $n=2$, then $R$ is an {\it {elliptic root system}} (see
\cite[(1.3) Example~1~iii)]{saito1},
\cite{saitotakebayashi} and
\cite{saitoyoshii}).

\begin{rem}\label{remark:MacSaito}
Assume $n=1$. Here 
we give a sketch of a proof of an equivalence between affine root systems in
the senses of 
\cite{MacDonald}, \cite[\S1.2]{MacDonald2} and \cite{saito1} (i.e. our Definition~\ref{definition:definitionexaff}). 
Let $F$ and $E$ be as in \cite[\S1.2]{MacDonald2}.
Let $S$ be a subset of $F$, and
assume $S$ is an irreducible affine root system in the sense of \cite[\S1.2]{MacDonald2}.
Identify $\v$ with $F$, that is, we regard $\v$ as
an $l+1$-dimensional $\bbbr$-linear space of affine-linear functions
$f:E\to\bbbr$. Clearly $S$ satisfies ({\rm{AX}}1) and ({\rm{AX}}3-5). 
Let $\lambda\in \vt$. Let $\mu
\in \vt$ be such that $c\mu\in \lambda+\vo$ for some $c\in\bbbrt$.
Then $\lambda-c\mu$ is a constant function on $E$,
that is, $(\lambda-c\mu)(E)=\{d_{\lambda-c\mu}\}$
for some $d_{\lambda-c\mu}\in\bbbr$.
We have $s_\mu s_\lambda(x)=x-(\lambda^\vee,x)(\lambda-c\mu)$
for $x\in \v$.
Further, for $e\in E$, we have 
$s_\mu s_\lambda\cdot e=e+{\frac {2d_{\lambda-c\mu}} {(\lambda,\lambda)}}D\lambda$,
see \cite[\S1.1]{MacDonald2} for $D\lambda$.
Then by using an argument similar to \cite[(1.16) Assertion~1]{saito1}, 
we can see that $S$
satisfies ({\rm{AX}}2).
Let $R$ be as in Definition~\ref{definition:definitionexaff}.
Let $T$ be
the subgroup of $W$ generated by $\{s_\al s_{\al^\prime}\,|\,
\al, \al^\prime\in R,\,\bbbrt\pi(\al)=\bbbrt\pi(\al^\prime)\}$.
Then $T$ 
is a normal abelian subgroup,
and $W/T$ can be identified with the finite Weyl group $W_{\pi(R)}$ (cf. \cite[(1.3)~Note 2 ii)]{saito1}). 
Then $R$ satisfies ({\rm{AR}}~4) of \cite[\S1.2]{MacDonald2}.
\end{rem}

\subsection{Base of an irreducible finite or affine root system}
   
Assume that $n\in\{0,1\}$ (cf. (\ref{eqn:ranknullitydef})). 
We call
a subset $\Pi$ of $R$ formed by $(l+n)$-linearly
independent elements a {\it base} if
\begin{equation} \label{eqn:eqndfofbase}
R=(R\cap\bbbzp \Pi)\cup(R\cap\bbbzm \Pi).
\end{equation}
(For $n=0$, see \cite[Theorem~10.1]{Humphreys}.
For $n=1$, see Theorem~\ref{theorem:theoremoneone}
in this paper 
(cf. MacDonald~\cite[(4.6)]{MacDonald} (see also
\cite[(3.3)~i)-iii)]{saito1})). If $\Pi$ is a base of $R$,
then, for $\bbbx \in\{\bbbz, \bbbr\}$, we have
\begin{equation}\label{eqn:eqnzrzpizal}
\mbox{$\Pi$ is an $\bbbx$-basis of $\bbbx R$, that is,}\,\,
\bbbx R=\bigoplus_{\al\in\Pi}\bbbx \al.
\end{equation} \par
Assume that $n=1$. Let
$\Pi=\{\al_0,\al_1,\ldots,\al_l\}$ be a base of $R$;
we always assume $\al_0$ is such that
$\{\pi(\al_1),\ldots,\pi(\al_l)\}$ is a
base of $\pi(R)$ (see
Theorem~\ref{theorem:theoremoneone}).
Let $\delta(\Pi)\in\bbbz \Pi$ be such that
\begin{equation}\label{eqn:eqndefondelpi}
\delta(\Pi)\in\bbbn \Pi\,\,\mbox{and $\{\delta(\Pi)\}$ is a $\bbbz$-basis of $(\bbbz R)^0$,
that is,}
\,\,\bbbz \delta(\Pi)=(\bbbz R)^0.
\end{equation}
$\delta(\Pi)$ is unique by (\ref{eqn:eqnbsofexaff}). 
By (\ref{eqn:bsofexaffd}), for $\bbbx\in\{\bbbz,\bbbr\}$,
we have
\begin{equation}\label{eqn:eqnzraff}
\mbox{$\{\al_1,\ldots,\al_l,\delta(\Pi)\}$ is a $\bbbx$-basis of $\bbbx R$, that is,}\,\,
\bbbx R
=(\bigoplus_{i=1}^n \bbbx \al_i)\oplus\bbbx \delta(\Pi).
\end{equation}
%Set
%\begin{equation}\label{eqn;eqndefofapi}
%A(\Pi):=((\al_i^\vee,\al_j))_{0\leq i,j\leq l}.
%\end{equation}
%Then $A(\Pi)$ is a generalized Cartan matrix of
%affine-type
%(see Theorem~\ref{theorem:theoremoneone}~(4)).
%If ${\widetilde \Pi}$ is another base of $R$,
%then $\delta({\widetilde \Pi})=\delta(\Pi)$ or
%$-\delta(\Pi)$
%and $A({\widetilde \Pi})=A(\Pi)$
%(after renumbering $\al_1,\ldots,\al_l$)
%(see Theorem~\ref{theorem:theoremoneone}~(3)-(4)).

The following lemma is well-known, e.g., see 
\cite[Theorem~10.3, Lemmas~10.4~C,D, \S12 Excercises~3]{Humphreys}.

\begin{lmm}\label{lemma:lemmatheta}
Assume that $n=0$ {\rm{(}}cf. {\rm{(\ref{eqn:ranknullitydef})}}{\rm{)}}.
Let $\Pi$ be a base of $R$ 
{\rm{(}}cf. {\rm{(\ref{eqn:eqndfofbase})}}{\rm{)}}. Then we have
the following{\rm{:}} \par
{\rm{(1)}} $W_\Pi=W$ and $W\cdot \Pi=R\setminus 2R$.
{\rm{(}}see {\rm{(\ref{eqn:eqndefofws})}} for $W_\Pi$
and see Definition~{\rm{\ref{definition:definitionexaff}}} for $W=W_R${\rm{)}}.\par
{\rm{(2)}} $W\cdot \al=\{\beta\in R |(\al,\al)=(\beta,\beta)\}$ for each $\al\in R$.
\par
{\rm{(3)}} For each $\al\in R$, there exists a unique 
$\al_+ \in W\cdot \al$ such that
$W\cdot \al\subset \al_++\bbbzm\Pi$.
\par
{\rm{(4)}} Let $r=|\{(\al,\al)|\al\in R\}|$. Then $1\leq r\leq 3$.
Moreover, if $r=3$, then $R\cap 2R=\{\beta\in R\,|\,(\beta,\beta)\geq (\al,\al)
\,\,\mbox{for all $\al\in R$}\}$.
\end{lmm}

{\it {Proof of {\rm{(3)}}}}. Let $\Pi=\{\al_1,\ldots,\al_l\}$. Then $\al_+$ is the element
$\sum_{i=1}^l m_i\al_i\in W\cdot \al$ ($m_i\in\bbbz$) for which $\sum_{i=1}^l m_i$ is maximal.
Let $w\in W_\Pi$ and let $w=s_{\al_1}\cdots s_{\al_r}$ be a reduced expression,
that is, $r$ is as small as possible. 
By \cite[Corollary~10.2~C]{Humphreys}, we have 
$w.\al_+=\al_+-\sum_{j=1}^r (\al_j^\vee,\al_+)s_{\al_1}\cdots s_{\al_{j-1}}(\al_j)
\in \al_++\bbbzm\Pi$. \hfill $\Box$
\newline\par
For $R$ and $\Pi$ of Lemma~\ref{lemma:lemmatheta},  we let
\begin{equation}\label{eqn:bigtheta}
\Theta(R,\Pi):=\{\al_+\in R|\al\in R\}.
\end{equation} By checking directly (and using \cite[\S12~Table~2]{Humphreys}),
we have
\begin{equation}\label{eqn:munumorezero} 
\mbox{$(\mu,\nu)>0$ for $\mu$, $\nu\in\Theta(R,\Pi)$}.
\end{equation} 

\subsection{Notation $S_{{\rm{sh}}}$, $S_{{\rm{lg}}}$, $S_{{\rm{ex}}}$}

Let $R$ be an ($n$-)extended affine root system 
(see Definition~\ref{definition:definitionexaff}).
Define the subsets $\Rsh$, $\Rlg$
and
$\Rex$ of $R$
by
\begin{equation}\nonumber
\Rsh:=\{\al\in R\,|\,(\al,\al)\leq (\beta,\beta)
\,\,\mbox{for all $\beta\in R$}\},
\end{equation} $\Rex:=R\cap \pi^{-1}(2
\pi(\Rsh))$
and $\Rlg :=R\setminus(\Rsh\cup\Rex)$
(see (\ref{eqn:eqnpipr}) for $\pi$). Then we have
\begin{equation}\label{eqn:rshrlgrex}
R=\Rsh \cup \Rlg \cup \Rex
\,\,\mbox{(disjoint union).}
\end{equation}
For a subset $S$ of $R$, let 
\begin{equation}\label{eqn:sshslgsex}
\mbox{$S_{{\rm{sh}}}:=S\cap
\Rsh $,
$S_{{\rm{lg}}}:=S\cap \Rlg$,
$S_{{\rm{ex}}}:=S\cap \Rex$.} 
\end{equation}

\section{A non-topological proof for the existence of a base of an affine root system}

In this section we assume $R$ is an affine root system, that is, we assume $n=1$
{\rm{(}}see  {\rm{(\ref{eqn:ranknullitydef}))}}.

\subsection{The existence of a base of an affine root system}

The following theorem seems to be well-known
(see \cite{MacDonald}), but we state and prove it for
later use. The proof in \cite{MacDonald} uses topological terminology.
Our proof seems to be the first one without using topology. 
Besides we need a technically written statement of the following theorem
for application.

\begin{thm}\label{theorem:theoremoneone}
{\rm{(}}cf. {\rm{\cite{MacDonald})}}
Let $\delta^\prime\in \vo\setminus\{0\}$ be such that $\bbbz \delta^\prime
=(\bbbz R)^0$
{\rm{(}}cf. {\rm{(\ref{eqn:eqnbsofexaff})}}{\rm{)}}. Let
$\Pi^\prime=\{\al_1,\ldots,\al_l\}$
be a subset of $R$ with $|\Pi^\prime|=l$ such that $\pi(\Pi^\prime)$ is
a base of the irreducible finite root system $(\pi(R), \v /\bbbr\delta^\prime)$
{\rm{(}}cf. {\rm{(\ref{eqn:eqndfofbase})}}
and {\rm{(\ref{eqn:eqnpipr})}}{\rm{)}}.
{\rm{(}}So $\bbbz R=\bbbz\delta^\prime\oplus\bbbz\Pi^\prime$
{\rm{(}}cf. {\rm{(\ref{eqn:bsofexaffd})}}{\rm{)}}.{\rm{)}}
Then
there exists a unique
\begin{equation}\label{eqn:defalo}
\al_0=\al_0(R,\Pi^\prime,\delta^\prime)\in R
\end{equation}
such that $\{\al_0\}\cup\Pi^\prime$ is a base of $R$ and
$\al_0\in\bbbn\delta^\prime\oplus\bbbz\Pi^\prime$.
Moreover $\al_0=\delta^\prime -\theta$ for some
$\theta\in\bbbn\Pi^\prime$ with
$\pi(\theta)\in\Theta(\pi(R),\pi(\Pi^\prime))$ 
{\rm{(}}see {\rm{(\ref{eqn:bigtheta})}}{\rm{)}}. 
In particular, $[(\al_i^\vee,\al_j)]_{0\leq i,j\leq l}$
is a generalized Cartan matrix
of affine-type in the sense of \cite[\S 4.3 and
Proposition~4.7]{kac1}.
Further, letting $\Pi_1=\{\al_0\}\cup \Pi^\prime$, for any base $\Pi_2$ of $R$ we have
$\Pi_2=\epsilon w(\Pi_1)$ for some
$\epsilon\in\{1,-1\}$ and $w\in W_{\Pi_1}$. In particular, 
\begin{equation}\label{eqn:wpipi}
R\setminus 2R=W_{\Pi_1}\cdot \Pi_1\,\,
\mbox{and}\,\,W=W_{\Pi_1}.
\end{equation}
\end{thm}
{\it{Proof.}}
({\it{Strategy.}} We use a linear map
$f:\v\rightarrow\bbbr$ (i.e., $f\in \v^*$)
such that $f(\al_i)=1$ ($1\leq i\leq l$) and
$f(\delta^\prime)$ is sufficiently
large (see (\ref{eqn:eqndefoff})). Let $\Pi^f$ be the
subset of $R$
formed by the elements $\beta\in R$
satisfying the condition that $f(\beta)>0$
and $\beta$ is not expressed as the summation of
more than one elements $\beta^\prime$ of $R$ with
$f(\beta^\prime)>0$ (see (\ref{eqn:eqndefpif})).
We show that $\Pi^f$ is a base of $R$ satisfying
the properties of the statement.
It is easy to see that $\Pi^\prime\subset\Pi^f$
and $R=(R\cap\bbbz_+\Pi^f)\cup(R\cap\bbbz_-\Pi^f))$.
We show $|\Pi^f|=l+1$ by using
(\ref{eqn:munumorezero}).)
\par
We proceed with the proof of the theorem in the following steps.
\par {\it{Step}}~1 ({\it{Definition of $f$}}).
Notice that for $\bbbx\in\{\bbbz,\bbbr\}$,
\begin{equation}\label{eqn:zdpopzai}
\bbbx R=\bbbx \delta^\prime
\oplus(\oplus_{i=1}^l\bbbx \al_i)
\end{equation} (see (\ref{eqn:bsofexaffd})).
We may assume that $(\al_i,\al_i)\leq (\al_{i+1},\al_{i+1})$
for $1\leq i\leq l-1$.
Also since $\pi(\Pi')$ is a base of $\pi(R)$, if $l\geq 2$, we
may assume $\al_1$ is such that
there exists a unique $j\in\{2,\ldots,l\}$
such that $(\al_1,\al_j)\ne 0$.
Let
\begin{equation}\label{eqn:defofrp}
R^\prime:=\left\{\begin{array}{ll}
W_{\Pi^\prime}\cdot (\Pi^\prime\cup\{2\al_1\}) &
\quad\mbox{if $l=1$,} \\
W_{\Pi^\prime}\cdot (\Pi^\prime\cup\{2\al_1\}) &
\quad\mbox{if $l\geq 2$ and $2(\al_1,\al_1)=(\al_2,\al_2)$,} \\
W_{\Pi^\prime}\cdot \Pi^\prime & \quad\mbox{otherwise.}
\end{array}\right.
\end{equation}
Using \cite[Theorem~10.3~(c)
(and \S12~Exercise~3)]{Humphreys},
we can see that 
$W_{\Pi^\prime}\cdot \Pi^\prime$ and $R^\prime$
are irreducible finite root systems
with the base $\Pi^\prime$.
If $\pi(R)$ is reduced, then $\pi(R)=\pi(W_{\Pi^\prime}\cdot \Pi^\prime)$.
If $\pi(R)$ is not reduced, then $\pi(R)=\pi(R^\prime)$.
In particular, we have
\begin{equation}\label{eqn:eqnrpzdp}
R\subset R^\prime +\bbbz \delta^\prime.
\end{equation} (see also (\ref{eqn:zdpopzai})).

Define $f\in \v ^*$ by
\begin{equation}\label{eqn:eqndefoff}
f(\al_i)=1\,\,(1\leq i\leq l)\quad\mbox{and}\quad
f(\delta^\prime)=3M,
\end{equation} where $M:=\max\{|f(\gamma)||\gamma\in
R^\prime\}$ (notice
$|R^\prime|<\infty$). It follows from (\ref{eqn:eqnrpzdp})
that
$f(\beta)\ne 0$ for $\beta\in R$.
\par {\it{Step}}~2 ({\it{Definition of $\Pi^f$}}).
Let $R^{f,+}:=\{\beta\in R|f(\beta)>0\}$.
By (\ref{eqn:eqndefoff}),
we have
\begin{equation}\label{eqn:eqncorfp}
R^{f,+}=R\cap((R^\prime\cap\bbbz _+\Pi^\prime)
\cup(\cup_{m=1}^\infty (m\delta^\prime+R^\prime))).
\end{equation}
Let
$\Pi^f$ be a subset of $R$ formed by the elements
$\beta\in R^{f,+}$ satisfying the condition that
there exist no $\beta_1,\ldots,\beta_r\in R^{f,+}$
with $r\geq 2$ such that $\beta=\beta_1+\cdots+\beta_r$;
namely,
\begin{equation}\label{eqn:eqndefpif}
\Pi^f:=R^{f,+}\setminus(\bigcup_{r=2}^\infty
\{\sum_{i=1}^r\beta_i|\beta_i\in R^{f,+}\}).
\end{equation}
By (\ref{eqn:eqncorfp}), we have
\begin{equation}\label{eqn:eqnppsubpf}
\Pi^\prime\subset\Pi^f.
\end{equation} Notice $\bbbz \Pi^\prime
\ne\bbbz R$ (by (\ref{eqn:zdpopzai})).
Then we have
\begin{equation}\label{eqn:eqnzrzpif}
\bbbz \Pi^f=\bbbz R,\,\,
R=(R\cap\bbbz _+\Pi^f)\cup (R\cap\bbbz _-\Pi^f)\,\,\mbox{and}\,\,
|\Pi^f|\geq |\Pi^\prime|+1.
\end{equation}
\par
(As mentioned in our strategy, we show that $\Pi^f$ is a
base of $R$.) \par
\par {\it{Step}}~3 ({\it{If $\beta\in\Pi^f/\Pi^\prime$,
then we have $\pi(\beta)\in\Theta(\pi(R),\pi(\Pi^\prime))$
{\rm{(}}for $\Theta(\pi(R),\pi(\Pi^\prime))$, see {\rm{(\ref{eqn:bigtheta})}}{\rm{)}}}}).
Let $\beta\in\Pi^f/\Pi^\prime$
(see also (\ref{eqn:eqnppsubpf})-(\ref{eqn:eqnzrzpif})).
We show that $\beta$ is expressed as
\begin{equation}\label{eqn:eqnmdmth}
\beta=m\delta^\prime-\theta
\end{equation} for some $m\in\bbbn$ and some
$\theta$ with
\begin{equation}\label{eqn:eqnttshlgex}
\theta\in\Theta(R^\prime,\Pi^\prime)
\end{equation} (see  (\ref{eqn:bigtheta})
for $\Theta(R^\prime,\Pi^\prime)$).
By (\ref{eqn:eqncorfp}),
since $\Pi^f\subset R^{f,+}$, we have
\begin{equation}\label{eqn:eqnmdpm}
\beta=m\delta^\prime+\mu
\end{equation} for some $m\in\bbbn$ and
$\mu\in R^\prime$.
Let $\theta\in\Theta(R^\prime,\Pi^\prime)
\cap W_{\Pi^\prime}\cdot \mu$,
where we recall from
Lemma~\ref{lemma:lemmatheta}~(2)-(3)
that $|\Theta(R^\prime,\Pi^\prime)
\cap W_{\Pi^\prime}\cdot \mu|=1$.
Notice $\{\mu,-\mu,\theta,-\theta\}\subset
W_{\Pi^\prime}\cdot \mu$
(cf. Lemma~\ref{lemma:lemmatheta}~(2)).
Then $m\delta^\prime-\theta\in R$ since
$m\delta^\prime-\theta\in
m\delta^\prime+W_{\Pi^\prime}\cdot \mu
=W_{\Pi^\prime}\cdot (m\delta^\prime+\mu)
=W_{\Pi^\prime}\cdot \beta\subset R$.
By
Lemma~\ref{lemma:lemmatheta}~(3),
we have $\theta+\mu=\theta-(-\mu)\in\bbbz _+\Pi^\prime$.
Since $m\delta^\prime-\theta\in R^{f,+}$ (cf.
(\ref{eqn:eqncorfp})),
$\beta=(m\delta^\prime-\theta)+(\theta+\mu)$
and $\beta\in\Pi^f$, we have $\theta+\mu=0$
and
(\ref{eqn:eqnmdmth}), as desired. \par
\par {\it{Step}}~4 ({\it{$|\Pi^f|=l+1$}}). We show
\begin{equation}\label{eqn:eqnpfsmppone}
|\Pi^f\setminus\Pi^\prime|=1,\,\,\mbox{i.e.,}\,\,|\Pi^f|=l+1
\end{equation} (see also (\ref{eqn:eqnppsubpf})-(\ref{eqn:eqnzrzpif})).
\par
Assume
$|\Pi^f\setminus\Pi^\prime|>1$.
Let $\beta_1$, $\beta_2\in\Pi^f\setminus\Pi^\prime$
and assume
$\beta_1\ne\beta_2$. Assume $(\beta_1,\beta_1)\leq
(\beta_2,\beta_2)$.
Then, by (\ref{eqn:munumorezero}) and
(\ref{eqn:eqnmdmth})-(\ref{eqn:eqnttshlgex}),
we see that
\begin{equation}\nonumber
(\beta_2^\vee,\beta_1)=
\left\{\begin{array}{ll}
1 & \mbox{if $\pi(\beta_1)\ne\pi(\beta_2)$,} \\
2 & \mbox{if $\pi(\beta_1)=\pi(\beta_2)$.}
\end{array}\right.
\end{equation}
Assume $(\beta_2^\vee,\beta_1)=1$. Then, since
$\pm(\beta_1-\beta_2)=s_{\beta_2}(\pm\beta_1)\in R$, we have
$\beta_1-\beta_2\in R^{f,+}$ or $\beta_2-\beta_1\in
R^{f,+}$.
This contradicts the fact $\beta_1$, $\beta_2\in\Pi^f$
since $\beta_1=\beta_2+(\beta_1-\beta_2)$
and $\beta_2=\beta_1+(\beta_2-\beta_1)$.
Assume $(\beta_2^\vee,\beta_1)=2$, so
$\pi(\beta_1)=\pi(\beta_2)$.
By (\ref{eqn:eqnmdmth}),
there exist $n_1$, $n_2\in\bbbn$ and $\theta
\in\Theta(R^\prime,\Pi^\prime)$
such that
\begin{equation}\nonumber
\beta_i=n_i\delta^\prime-\theta
\quad(i\in\{1,2\})
\end{equation} (so $\beta_2-\beta_1=(n_2-n_1)\delta^\prime$).
Assume
$n_1<n_2$. Notice that for $i\in\{1,2\}$
and $r\in\bbbz $,
\begin{eqnarray}\label{eqn:eqnrnisbonesbtwopre}
R&\ni& (s_{\beta_2}s_{\beta_1})^r(\beta_i) \quad (\mbox{by ({\rm{AX}}4)})\\
&=& (n_i+2r(n_2-n_1))\delta^\prime - \theta \nonumber\\
&=&
\left\{\begin{array}{ll}
(n_2+(2r-1)(n_2-n_1))\delta^\prime - \theta &
\mbox{if $i=1$,} \nonumber\\
(n_2+2r(n_2-n_1))\delta^\prime - \theta &
\mbox{if $i=2$.}
\end{array}\right. \nonumber
\end{eqnarray} Hence
\begin{eqnarray}\label{eqn:eqnrnisbonesbtwo}
(n_2+r(n_2-n_1))\delta^\prime-\theta\in R
\quad\mbox{for all $r\in\bbbz $.}
\end{eqnarray}
Let
$n_3\in\bbbz _+$ and $t\in\bbbn$ be such that
$0\leq n_3<n_2-n_1$ and
$n_2=t(n_2-n_1)+n_3$.
Assume $n_3=0$.
By (\ref{eqn:eqnrnisbonesbtwo}),
$\{-\theta,(n_2-n_1)\delta^\prime-\theta\}\subset R$.
Hence, by
(\ref{eqn:eqncorfp})
(and (\ref{eqn:eqnoneone})),
$\{\theta,(n_2-n_1)\delta^\prime-\theta\}\subset
R^{f,+}$.
Notice $t\geq 2$ (since $0<n_1<n_2$ and $n_3=0$). Since
$\beta_2=t((n_2-n_1)\delta^\prime-\theta)+(t-1)\theta$, we have
$\beta_2\notin \Pi^f$, contradiction.
Assume $n_3>0$.
Notice $2n_3< n_2$
(since $2n_3<(n_2-n_1)+n_3\leq t(n_2-n_1)+n_3
=n_2$). Let $\beta_3=n_3\delta^\prime-\theta$.
By (\ref{eqn:eqnrnisbonesbtwo}), $\beta_3\in R$.
By (\ref{eqn:eqncorfp}), $\beta_3\in R^{f,+}$.
Notice $\beta_2-2\beta_3=s_{\beta_3}(\beta_2)
\in R$ (by ({\rm{AX}}4)). Then
by (\ref{eqn:eqncorfp}), we have
\begin{equation}\nonumber
\beta_2-2\beta_3=(n_2-2n_3)\delta^\prime+\theta\in R^{f,+}.
\end{equation} Since
$\beta_2=(\beta_2-2\beta_3)+2\beta_3$, we have
$\beta_2\notin \Pi^f$, contradiction. Hence
$|\Pi^f|=l+1$, as desired.
\par
\par {\it{Step}}~5 ({\it{$\Pi^f$ is a base with $\al_0=\delta^\prime-\theta$}}).
Let $\al_0$ be $\beta=m\delta^\prime-\theta$ of (\ref{eqn:eqnmdmth}).
Then $\Pi^f=\Pi^\prime\cup\{\al_0\}$,
where we notice (\ref{eqn:eqnppsubpf}) and (\ref{eqn:eqnpfsmppone}).
It is clear that the elements of $\Pi^f$
are linearly independent (cf. (\ref{eqn:zdpopzai})).
Hence, by (\ref{eqn:eqnzrzpif}),
$\Pi^f$ is a base of $R$
(cf. (\ref{eqn:eqndfofbase})).
Since $\bbbz\Pi^\prime\oplus\bbbz\delta^\prime
=\bbbz\Pi^\prime\oplus\bbbz\al_0$
(by (\ref{eqn:zdpopzai}) and (\ref{eqn:eqnzrzpif})),
we have $m=1$. 

\par {\it{Step}}~6 ({\it{The last claim holds}}).
Let $\Pi_1=\Pi^\prime\cup\{\al_0\}$.
Let $\Pi_2$ be a base of $R$.
Define $h\in V^*$ by $h(\beta):=1$ ($\beta\in\Pi_2$).
Then $h(R)\subset\bbbz\setminus\{0\}$.
By the same formula as in
(\ref{eqn:eqnrnisbonesbtwopre}),
we have $|\{(s_\theta s_{\al_0})^r(\al_0)\in
R|r\in\bbbz\}|
=\infty$ (notice that $(s_\theta s_{\al_0})^r(\al_0)\in
R$ (by ({\rm{AX}}4))
since $s_\theta=s_{{\frac 1 2}\theta}$ and $\theta\in
R\cup 2R$ (see (\ref{eqn:eqnttshlgex})
and (\ref{eqn:defofrp}))).
Hence $|R|=\infty$, which implies $|h(R)|=\infty$.
Hence, by (\ref{eqn:eqnrpzdp}), since $|R^\prime|<\infty$
($R^\prime$ is an irreducible finite root system),
we have $h(\delta^\prime)\ne 0$. We may
assume
\begin{equation}\label{eqn:eqnrdpos}
h(\delta^\prime)>0
\end{equation}
(otherwise, we replace $\Pi_2$ with $-\Pi_2$).
Let
\begin{eqnarray*}\nonumber
m(\Pi_1,\Pi_2)&:=&|(R\cap\bbbz_+\Pi_1\cap\bbbz_-\Pi_2) 
\setminus 2R| \\
&=& |\{\beta\in
(R\cap\bbbz_+\Pi_1)\setminus
2R\,|\,h(\beta)<0\}|.
\end{eqnarray*} Since $\al_0=\delta^\prime-\theta$,
we have
$R\cap\bbbz_+\Pi_1\subset R^\prime+\bbbz_+\delta^\prime$
(cf. (\ref{eqn:eqnrpzdp})). Hence, since $|R^\prime|<\infty$,
by (\ref{eqn:eqnrdpos}),
we have $m(\Pi_1,\Pi_2)<\infty$.
%\begin{equation}\nonumber
%m(\Pi_1,\Pi_2)=|\{\beta\in
%(R\cap\bbbz_+\Pi_1)\setminus
%2R\,|\,h(\beta)<0\}|<\infty.
%\end{equation} 
\par
We use induction on $m(\Pi_1,\Pi_2)$; if
$m(\Pi_1,\Pi_2)=0$, then,
by (\ref{eqn:eqndfofbase}),  
$R\cap\bbbz_+\Pi_1=R\cap\bbbz_+\Pi_2$,
so $\Pi_1=\Pi_2$.
Assume $m(\Pi_1,\Pi_2)>0$.
Then there exists $\al\in\Pi_1$ such that
$\al\in\bbbz_-\Pi_2$
(notice that $R\subset\bbbz_-\Pi_2\cup\bbbz_+\Pi_2$).
By (\ref{eqn:eqndfofbase}) (and (\ref{eqn:eqnoneone})), we see
\begin{equation}\label{eqn:eqnsalrcapzpp}
s_\al((R\cap\bbbz_+\Pi_1)\setminus 2R)
=\{-\al\}\cup(((R\cap\bbbz_+\Pi_1)\setminus
2R)\setminus\{\al\}).
\end{equation} Then we have
\begin{eqnarray*}
\lefteqn{m(\Pi_1,s_\al(\Pi_2))} \\
&=&|(R\cap\bbbz_+\Pi_1\cap
\bbbz_-s_\al(\Pi_2))
\setminus 2R| \\
&=&|s_\al((R\cap\bbbz_+\Pi_1\cap
\bbbz_-s_\al(\Pi_2))
\setminus 2R)| \\
&=&|(s_\al(R\cap\bbbz_+\Pi_1)\cap
\bbbz_-\Pi_2)
\setminus 2R| \\
&=& m(\Pi_1,\Pi_2)-1
\quad\mbox{(by (\ref{eqn:eqnsalrcapzpp})
since $s_\al(\al)=-\al\notin\bbbz_-\Pi_2$)}.
\end{eqnarray*} Then, by the induction,
we see that there exists $w\in W_{\Pi_1}$ such that
$w(\Pi_2)=\Pi_1$, as desired.

Note that for any $\beta\in R\setminus 2R$,
there exists a subset $\Pi^{\prime\prime}$ of $R$
with $|\Pi^{\prime\prime}|=l$
such that $\beta\in\Pi^{\prime\prime}$
and $\pi(\Pi^{\prime\prime})$ is a base of $\pi(R)$.
Hence by the above argument, we have (\ref{eqn:wpipi}). 
This completes the proof. \hfill $\Box$
\newline\par
By (\ref{eqn:wpipi}), we have 
\begin{equation}\label{eqn:affnonrt} 
\left\{\begin{array}{l} 
R= W_\Pi\cdot (\Pi\cup(2\Pi\cap R)),\\
(\bbbz R)^\times\setminus R \\ \quad
=W_\Pi\cdot \Bigl((2\Pi\setminus R)\cup(\bigcup_{r\in 3+\bbbzp}r\Pi)\cup
((\bbbz R)^\times\setminus(\bbbzp\Pi\cup\bbbzm\Pi))\Bigr).
\end{array}\right.
\end{equation}

\subsection{Dynkin diagrams of affine root systems}
\label{subsection:superaffdynkin} 
Here we give the Dynkin diagrams for $(R,\Pi)$ of Theorem~\ref{theorem:theoremoneone}.
We assume that if $2\al_0\in R$, then $2\al_i\in R$ for some $i\ne 0$,
see $A^{(4)}(0,2l)$ below.
We describe them in 
the same manner as in \cite[Table 1-4]{kac2};
especially, if $2\al_i\notin R$ (resp. $2\al_i\in R$ ),
then the $i$-th dot is white (resp. black).
The names of them are also the same as in \cite[Table 1-4]{kac2}.
\newline\par
(i) The case of $l=1$:
\begin{center}
$A_1^{(1)}$
\setlength{\unitlength}{1mm}
\begin{picture}(18,8)(0,0)

\put(1,4.5){$\al_1$}\put(11.5,4.5){$\al_0$}

\put(2, 1){\circle{3}}
\put(12.5,1){\circle{3}}
\put(3.5, 0){$\Longleftrightarrow$}

\end{picture} $A_2^{(2)}$
\setlength{\unitlength}{1mm}
\begin{picture}(18,8)(0,0)

\put(1,4.5){$\al_1$}\put(11.5,4.5){$\al_0$}

\put(2, 1){\circle{3}}
\put(12,1){\circle{3}}
\put(4.5,-0.5){\line(1,0){5.5}}
\put(4.5,0.5){\line(1,0){5.5}}
\put(4.5,1.5){\line(1,0){5.5}}
\put(4.5,2.5){\line(1,0){5.5}}
\put(3.5,0){$\langle$}

\end{picture}

\end{center}\begin{center}

$B^{(1)}(0,1)$
\setlength{\unitlength}{1mm}
\begin{picture}(18,8)(0,0)
\put(1,4.5){$\al_1$}\put(11.5,4.5){$\al_0$}

\put(2, 1){\circle*{3}}
\put(12,1){\circle{3}}
\put(4.5,-0.5){\line(1,0){5.5}}
\put(4.5,0.5){\line(1,0){5.5}}
\put(4.5,1.5){\line(1,0){5.5}}
\put(4.5,2.5){\line(1,0){5.5}}
\put(3.5,0){$\langle$}

\end{picture} $C^{(2)}(2)$
\setlength{\unitlength}{1mm}
\begin{picture}(18,8)(0,0)
\put(1,4.5){$\al_1$}\put(11.5,4.5){$\al_0$}

\put(2, 1){\circle*{3}}
\put(12.5,1){\circle*{3}}
\put(3.5, 0){$\Longleftrightarrow$}

\end{picture} $A^{(4)}(0,2)$
\setlength{\unitlength}{1mm}
\begin{picture}(18,8)(0,0)
\put(1,4.5){$\al_1$}\put(11.5,4.5){$\al_0$}

\put(2, 1){\circle*{3}}
\put(12.5,1){\circle{3}}
\put(3.5, 0){$\Longleftrightarrow$}

\end{picture}

\end{center}
\vspace{1cm}

(ii) The case of $l=2$:
\begin{center}
$A_2^{(1)}$
\setlength{\unitlength}{1mm}
\begin{picture}(21,8)(0,0)
\put(-1,4.5){$\al_1$}\put(13,4.5){$\al_2$}
\put(6.5,10.5){$\al_0$}

\put(2, 1){\circle{3}}\put(8, 7){\circle{3}}
\put(14,1){\circle{3}}
\put(3.5, 1){\line(1,0){9}}
\put(3, 2){\line(1,1){4}}\put(13, 2){\line(-1,1){4}}

\end{picture}
$C_2^{(1)}$
\setlength{\unitlength}{1mm}
\begin{picture}(28,8)(0,0)
\put(1,4.5){$\al_2$}\put(11,4.5){$\al_1$}
\put(21,4.5){$\al_0$}

\put(2, 1){\circle{3}}
\put(12,1){\circle{3}}
\put(3.5, 0){$\Longrightarrow$}
\put(22,1){\circle{3}}
\put(13.5, 0){$\Longleftarrow$}

\end{picture}
$G_2^{(1)}$
\setlength{\unitlength}{1mm}
\begin{picture}(28,18)(0,0)
\put(1,4.5){$\al_1$}\put(11.5,4.5){$\al_2$}
\put(21.5,4.5){$\al_0$}

\put(2, 1){\circle{3}}
\put(12,1){\circle{3}}\put(22,1){\circle{3}}
\put(4.5,0){\line(1,0){5.5}}
\put(4.5,1){\line(1,0){5.5}}
\put(4.5,2){\line(1,0){5.5}}
\put(3.5,0){$\langle$}
\put(13.5,1){\line(1,0){7}}
\end{picture}

\end{center}
\begin{center}

$A_4^{(2)}$
\setlength{\unitlength}{1mm}
\begin{picture}(28,8)(0,0)
\put(1,4.5){$\al_1$}\put(11,4.5){$\al_2$}
\put(21,4.5){$\al_0$}

\put(2, 1){\circle{3}}
\put(12,1){\circle{3}}
\put(3.5, 0){$\Longleftarrow$}
\put(22,1){\circle{3}}
\put(13.5, 0){$\Longleftarrow$}

\end{picture}
$D_3^{(2)}$
\setlength{\unitlength}{1mm}
\begin{picture}(28,8)(0,0)
\put(1,4.5){$\al_1$}\put(11,4.5){$\al_2$}
\put(21,4.5){$\al_0$}

\put(2, 1){\circle{3}}
\put(12,1){\circle{3}}
\put(3.5, 0){$\Longleftarrow$}
\put(22,1){\circle{3}}
\put(13.5, 0){$\Longrightarrow$}

\end{picture}$D_4^{(3)}$
\setlength{\unitlength}{1mm}
\begin{picture}(28,8)(0,0)
\put(1,4.5){$\al_0$}\put(11.5,4.5){$\al_1$}
\put(21.5,4.5){$\al_2$}

\put(2, 1){\circle{3}}
\put(12,1){\circle{3}}\put(22,1){\circle{3}}
\put(3.5,1){\line(1,0){7}}
\put(14.5,0){\line(1,0){5.5}}
\put(14.5,1){\line(1,0){5.5}}
\put(14.5,2){\line(1,0){5.5}}
\put(13.5,0){$\langle$}

\end{picture}
\end{center}
\begin{center}

$B^{(1)}(0,2)$
\setlength{\unitlength}{1mm}
\begin{picture}(28,8)(0,0)
\put(1,4.5){$\al_1$}\put(11,4.5){$\al_2$}
\put(21,4.5){$\al_0$}

\put(2, 1){\circle*{3}}
\put(12,1){\circle{3}}
\put(3.5, 0){$\Longleftarrow$}
\put(22,1){\circle{3}}
\put(13.5, 0){$\Longleftarrow$}

\end{picture}$A^{(2)}(0,3)$
\setlength{\unitlength}{1mm}
\begin{picture}(28,8)(0,0)
\put(1,4.5){$\al_2$}\put(11,4.5){$\al_1$}
\put(21,4.5){$\al_0$}

\put(2, 1){\circle{3}}
\put(12,1){\circle*{3}}
\put(3.5, 0){$\Longrightarrow$}
\put(22,1){\circle{3}}
\put(13.5, 0){$\Longleftarrow$}

\end{picture}\end{center}\begin{center}
$C^{(2)}(3)$
\setlength{\unitlength}{1mm}
\begin{picture}(28,8)(0,0)
\put(1,4.5){$\al_1$}\put(11,4.5){$\al_2$}
\put(21,4.5){$\al_0$}

\put(2, 1){\circle*{3}}
\put(12,1){\circle{3}}
\put(3.5, 0){$\Longleftarrow$}
\put(22,1){\circle*{3}}
\put(13.5, 0){$\Longrightarrow$}

\end{picture}$A^{(4)}(0,4)$
\setlength{\unitlength}{1mm}
\begin{picture}(28,8)(0,0)
\put(1,4.5){$\al_1$}\put(11,4.5){$\al_2$}
\put(21,4.5){$\al_0$}

\put(2, 1){\circle*{3}}
\put(12,1){\circle{3}}
\put(3.5, 0){$\Longleftarrow$}
\put(22,1){\circle{3}}
\put(13.5, 0){$\Longrightarrow$}

\end{picture}
\end{center}
\vspace{1cm}

(iii) The case of $l\geq 3$:

\begin{center}
$D_{l+1}^{(2)}$
\setlength{\unitlength}{1mm}
\begin{picture}(57,7)(0,4)

\put(00.5, 0){$\al_1$}
\put(10.5, 0){$\al_2$}
\put(20.5, 0){$\al_3$}
\put(42.5, 0){$\al_l$}
\put(52.5, 0){$\al_0$}

\put(02,5){\circle{3}}
\put(03.5,4){$\Longleftarrow$}
\put(12,5){\circle{3}}
\put(13.5, 5){\line(1,0){7}}
\put(22,5){\circle{3}}
\put(23.5, 5){\line(1,0){5}}
\put(30.5, 4){$\cdots$}

\put(37.5, 5){\line(1,0){5}}
\put(44,5){\circle{3}}
\put(45.5,4){$\Longrightarrow$}

\put(54,5){\circle{3}}

\end{picture}\end{center}
\begin{center}
$C^{(2)}(l+1)$
\setlength{\unitlength}{1mm}
\begin{picture}(57,7)(0,4)

\put(00.5, 0){$\al_1$}
\put(10.5, 0){$\al_2$}
\put(20.5, 0){$\al_3$}
\put(42.5, 0){$\al_l$}
\put(52.5, 0){$\al_0$}

\put(02,5){\circle*{3}}
\put(03.5,4){$\Longleftarrow$}
\put(12,5){\circle{3}}
\put(13.5, 5){\line(1,0){7}}
\put(22,5){\circle{3}}
\put(23.5, 5){\line(1,0){5}}
\put(30.5, 4){$\cdots$}

\put(37.5, 5){\line(1,0){5}}
\put(44,5){\circle{3}}
\put(45.5,4){$\Longrightarrow$}

\put(54,5){\circle*{3}}

\end{picture}
\end{center}
\begin{center}
$A^{(4)}(0,2l)$
\setlength{\unitlength}{1mm}
\begin{picture}(57,7)(0,4)

\put(00.5, 0){$\al_1$}
\put(10.5, 0){$\al_2$}
\put(20.5, 0){$\al_3$}
\put(42.5, 0){$\al_l$}
\put(52.5, 0){$\al_0$}

\put(02,5){\circle*{3}}
\put(03.5,4){$\Longleftarrow$}
\put(12,5){\circle{3}}
\put(13.5, 5){\line(1,0){7}}
\put(22,5){\circle{3}}
\put(23.5, 5){\line(1,0){5}}
\put(30.5, 4){$\cdots$}

\put(37.5, 5){\line(1,0){5}}
\put(44,5){\circle{3}}
\put(45.5,4){$\Longrightarrow$}

\put(54,5){\circle{3}}

\end{picture}
\end{center}
\begin{center}
$A_{2l}^{(2)}$
\setlength{\unitlength}{1mm}
\begin{picture}(57,7)(0,4)

\put(00.5, 0){$\al_1$}
\put(10.5, 0){$\al_2$}
\put(20.5, 0){$\al_3$}
\put(42.5, 0){$\al_l$}
\put(52.5, 0){$\al_0$}

\put(02,5){\circle{3}}
\put(03.5,4){$\Longleftarrow$}
\put(12,5){\circle{3}}
\put(13.5, 5){\line(1,0){7}}
\put(22,5){\circle{3}}
\put(23.5, 5){\line(1,0){5}}
\put(30.5, 4){$\cdots$}

\put(37.5, 5){\line(1,0){5}}
\put(44,5){\circle{3}}
\put(45.5,4){$\Longleftarrow$}

\put(54,5){\circle{3}}

\end{picture} \end{center}
\begin{center} $B^{(1)}(0,l)$
\setlength{\unitlength}{1mm}
\begin{picture}(57,7)(0,4)

\put(00.5, 0){$\al_1$}
\put(10.5, 0){$\al_2$}
\put(20.5, 0){$\al_3$}
\put(42.5, 0){$\al_l$}
\put(52.5, 0){$\al_0$}

\put(02,5){\circle*{3}}
\put(03.5,4){$\Longleftarrow$}
\put(12,5){\circle{3}}
\put(13.5, 5){\line(1,0){7}}
\put(22,5){\circle{3}}
\put(23.5, 5){\line(1,0){5}}
\put(30.5, 4){$\cdots$}

\put(37.5, 5){\line(1,0){5}}
\put(44,5){\circle{3}}
\put(45.5,4){$\Longleftarrow$}

\put(54,5){\circle{3}}

\end{picture}
\end{center}

\begin{center} 

$B_l^{(1)}$
\setlength{\unitlength}{1mm}
\begin{picture}(47,17)(0,4)

\put(36, 13){$\al_0$}
\put(00.5, 0){$\al_1$}
\put(10.5, 0){$\al_2$}
\put(32.5, 0){$\al_{l-1}$}
\put(42.5, 0){$\al_l$}

\put(02,5){\circle{3}}
\put(03.5,4){$\Longleftarrow$}
\put(12,5){\circle{3}}
\put(34,6.5){\line(0,1){7}}\put(34,15){\circle{3}}
\put(13.5, 5){\line(1,0){5}}
\put(20.5, 4){$\cdots$}

\put(27.5, 5){\line(1,0){5}}
\put(34,5){\circle{3}}
\put(35.5,5){\line(1,0){7}}

\put(44,5){\circle{3}}

\end{picture} \par
$A^{(2)}(0,2l-1)$
\setlength{\unitlength}{1mm}
\begin{picture}(47,17)(0,4)

\put(36, 13){$\al_0$}
\put(00.5, 0){$\al_1$}
\put(10.5, 0){$\al_2$}
\put(32.5, 0){$\al_{l-1}$}
\put(42.5, 0){$\al_l$}

\put(02,5){\circle*{3}}
\put(03.5,4){$\Longleftarrow$}
\put(12,5){\circle{3}}
\put(34,6.5){\line(0,1){7}}\put(34,15){\circle{3}}
\put(13.5, 5){\line(1,0){5}}
\put(20.5, 4){$\cdots$}

\put(27.5, 5){\line(1,0){5}}
\put(34,5){\circle{3}}
\put(35.5,5){\line(1,0){7}}

\put(44,5){\circle{3}}

\end{picture}
\end{center}

\begin{center} $A_l^{(1)}$
\setlength{\unitlength}{1mm}
\begin{picture}(47,17)(0,4)

\put(24.5, 15){$\al_0$}
\put(00.5, 0){$\al_1$}
\put(10.5, 0){$\al_2$}

\put(30.5, 0){$\al_{l-1}$}
\put(40.5, 0){$\al_l$}

\put(03.2,6.2){\line(2,1){17}}
\put(02,5){\circle{3}}
\put(03.5,5){\line(1,0){7}}
\put(12,5){\circle{3}}

\put(13.5, 5){\line(1,0){4}}\put(19.5, 4){$\cdots$}
\put(22,15){\circle{3}}
\put(26.5,5){\line(1,0){4}}
\put(32,5){\circle{3}}
\put(33.5,5){\line(1,0){7}}
\put(42,5){\circle{3}}\put(40.8,6.2){\line(-2,1){17}}

\end{picture}
$D_l^{(1)}$
\setlength{\unitlength}{1mm}
\begin{picture}(47,17)(0,4)

\put(14, 13){$\al_0$}
\put(36, 13){$\al_l$}
\put(00.5, 0){$\al_1$}
\put(10.5, 0){$\al_2$}
\put(32.5, 0){$\al_{l-2}$}
\put(42.5, 0){$\al_{l-1}$}

\put(12,6.5){\line(0,1){7}}\put(12,15){\circle{3}}
\put(02,5){\circle{3}}
\put(03.5,5){\line(1,0){7}}
\put(12,5){\circle{3}}
\put(34,6.5){\line(0,1){7}}\put(34,15){\circle{3}}
\put(13.5, 5){\line(1,0){5}}
\put(20.5, 4){$\cdots$}

\put(27.5, 5){\line(1,0){5}}
\put(34,5){\circle{3}}
\put(35.5,5){\line(1,0){7}}

\put(44,5){\circle{3}}

\end{picture}
\end{center}
\begin{center}
$C_l^{(1)}$
\setlength{\unitlength}{1mm}
\begin{picture}(57,7)(0,4)

\put(00.5, 0){$\al_1$}
\put(10.5, 0){$\al_2$}
\put(20.5, 0){$\al_3$}
\put(42.5, 0){$\al_l$}
\put(52.5, 0){$\al_0$}

\put(02,5){\circle{3}}
\put(03.5,4){$\Longrightarrow$}
\put(12,5){\circle{3}}
\put(13.5, 5){\line(1,0){7}}
\put(22,5){\circle{3}}
\put(23.5, 5){\line(1,0){5}}
\put(30.5, 4){$\cdots$}

\put(37.5, 5){\line(1,0){5}}
\put(44,5){\circle{3}}
\put(45.5,4){$\Longleftarrow$}

\put(54,5){\circle{3}}

\end{picture} 
\end{center}
\begin{center} $E_6^{(1)}$
\setlength{\unitlength}{1mm}
\begin{picture}(47,27)(0,4)

\put(24, 13){$\al_6$}\put(24, 23){$\al_0$}
\put(00.5, 0){$\al_1$}
\put(10.5, 0){$\al_2$}
\put(20.5, 0){$\al_3$}
\put(30.5, 0){$\al_4$}
\put(40.5, 0){$\al_5$}

\put(02,5){\circle{3}}
\put(03.5,5){\line(1,0){7}}
\put(12,5){\circle{3}}
\put(22,6.5){\line(0,1){7}}\put(22,15){\circle{3}}
\put(13.5, 5){\line(1,0){7}}
\put(22,5){\circle{3}}
\put(23.5,5){\line(1,0){7}}
\put(32,5){\circle{3}}
\put(33.5,5){\line(1,0){7}}
\put(42,5){\circle{3}}
\put(22,16.5){\line(0,1){7}}\put(22,25){\circle{3}}
\end{picture}$A_{2l-1}^{(2)}$
\setlength{\unitlength}{1mm}
\begin{picture}(47,17)(0,4)

\put(36, 13){$\al_0$}
\put(00.5, 0){$\al_l$}
\put(10.5, 0){$\al_{l-1}$}
\put(32.5, 0){$\al_2$}
\put(42.5, 0){$\al_1$}

\put(02,5){\circle{3}}
\put(03.5,4){$\Longrightarrow$}
\put(12,5){\circle{3}}
\put(34,6.5){\line(0,1){7}}\put(34,15){\circle{3}}
\put(13.5, 5){\line(1,0){5}}
\put(20.5, 4){$\cdots$}

\put(27.5, 5){\line(1,0){5}}
\put(34,5){\circle{3}}
\put(35.5,5){\line(1,0){7}}

\put(44,5){\circle{3}}

\end{picture}
\end{center}

\begin{center} $E_7^{(1)}$
\setlength{\unitlength}{1mm}
\begin{picture}(67,17)(0,4)

\put(34, 13){$\al_7$}
\put(0.5, 0){$\al_0$}
\put(10.5, 0){$\al_1$}
\put(20.5, 0){$\al_2$}
\put(30.5, 0){$\al_3$}
\put(40.5, 0){$\al_4$}
\put(50.5, 0){$\al_5$}
\put(60.5, 0){$\al_6$}

\put(02,5){\circle{3}}
\put(03.5,5){\line(1,0){7}}
\put(12,5){\circle{3}}
\put(13.5,5){\line(1,0){7}}
\put(22,5){\circle{3}}
\put(32,6.5){\line(0,1){7}}\put(32,15){\circle{3}}
\put(23.5, 5){\line(1,0){7}}
\put(32,5){\circle{3}}
\put(33.5,5){\line(1,0){7}}
\put(42,5){\circle{3}}
\put(43.5,5){\line(1,0){7}}
\put(52,5){\circle{3}}
\put(53.5,5){\line(1,0){7}}
\put(62,5){\circle{3}}

\end{picture}
\end{center}

\begin{center} $E_8^{(1)}$
\setlength{\unitlength}{1mm}
\begin{picture}(77,17)(0,4)

\put(54, 13){$\al_8$}
\put(0.5, 0){$\al_0$}
\put(10.5, 0){$\al_1$}
\put(20.5, 0){$\al_2$}
\put(30.5, 0){$\al_3$}
\put(40.5, 0){$\al_4$}
\put(50.5, 0){$\al_5$}
\put(60.5, 0){$\al_6$}
\put(70.5, 0){$\al_7$}

\put(02,5){\circle{3}}
\put(03.5,5){\line(1,0){7}}
\put(12,5){\circle{3}}
\put(13.5,5){\line(1,0){7}}
\put(22,5){\circle{3}}
\put(52,6.5){\line(0,1){7}}\put(52,15){\circle{3}}
\put(23.5, 5){\line(1,0){7}}
\put(32,5){\circle{3}}
\put(33.5,5){\line(1,0){7}}
\put(42,5){\circle{3}}
\put(43.5,5){\line(1,0){7}}
\put(52,5){\circle{3}}
\put(53.5,5){\line(1,0){7}}
\put(62,5){\circle{3}}
\put(63.5,5){\line(1,0){7}}
\put(72,5){\circle{3}}

\end{picture}
\end{center}

\begin{center} $F_4^{(1)}$
\setlength{\unitlength}{1mm}
\begin{picture}(47,7)(0,4)

\put(00.5, 0){$\al_0$}
\put(10.5, 0){$\al_4$}
\put(20.5, 0){$\al_3$}
\put(30.5, 0){$\al_2$}
\put(40.5, 0){$\al_1$}

\put(02,5){\circle{3}}
\put(03.5,5){\line(1,0){7}}
\put(12,5){\circle{3}}

\put(13.5, 5){\line(1,0){7}}
\put(22,5){\circle{3}}
\put(23.5,4){$\Longrightarrow$}
\put(32,5){\circle{3}}
\put(33.5,5){\line(1,0){7}}
\put(42,5){\circle{3}}

\end{picture}
    $E_6^{(2)}$
\setlength{\unitlength}{1mm}
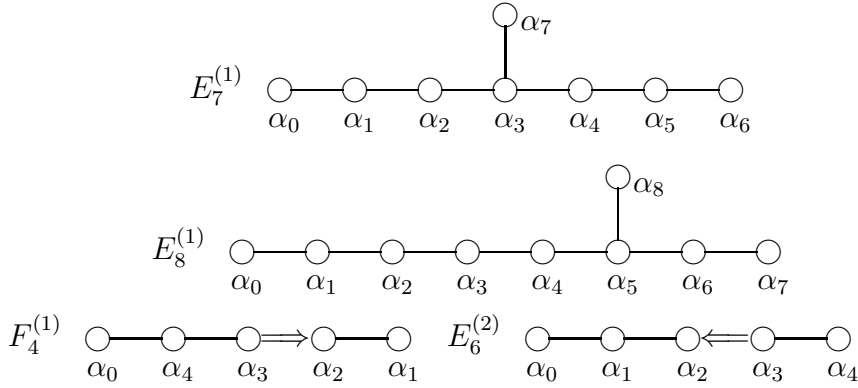
\begin{picture}(47,7)(0,4)

\put(00.5, 0){$\al_0$}
\put(10.5, 0){$\al_1$}
\put(20.5, 0){$\al_2$}
\put(30.5, 0){$\al_3$}
\put(40.5, 0){$\al_4$}

\put(02,5){\circle{3}}
\put(03.5,5){\line(1,0){7}}
\put(12,5){\circle{3}}

\put(13.5, 5){\line(1,0){7}}
\put(22,5){\circle{3}}
\put(23.5,4){$\Longleftarrow$}
\put(32,5){\circle{3}}
\put(33.5,5){\line(1,0){7}}
\put(42,5){\circle{3}}

\end{picture}\newline\newline
\end{center}

\section{Elliptic root systems}\label{section:elliptic}

In this section we assume $R$ is a reduced elliptic root system, that is, 
$R\cap 2R=\emptyset$ and 
$n=2$
{\rm{(}}see  {\rm{(\ref{eqn:ranknullitydef}))}}.

\subsection{Fundamental-set of an elliptic root system}\label{subsection:fundamental}
\begin{dfn}\label{definition:defoffset}
{\rm{({\it{Fundamental-set
$\Pi\cup\{a\}$}})
We say that a subset $\Pi\cup\{a\}$ of $\bbbz R$ is 
a {\it{fundamental-set}} of $R$ if
it satisfies the axioms
({\rm{FS}}1)-({\rm{FS}}2) below;
we always let
\begin{equation}\label{eqn:eqnpia}
\pi_a:\v \rightarrow \v/\bbbr a
\end{equation} denote the canonical map.
\newline\par
({\rm{FS}}1) $a\in(\bbbz R)^0$
and there exists
$b\in(\bbbz R)^0 $ such that
$\{a,b\}$ is a basis of $(\bbbz R)^0$, i.e.,
$(\bbbz R)^0 =\bbbz a\oplus\bbbz b$.
\par
({\rm{FS}}2) $|\Pi|=l+1$,
$\Pi\subset R$ and
$\pi_a(\Pi)$ is a base of the affine root
system $\pi_a(R)$.
}}
\end{dfn} \par Until end of this section, let $\Pi\cup\{a\}=\{\al_0,\ldots,\al_l\}\cup\{a\}$ denote a fundamental-set of $R$. We assume $\pi(\{\al_1,\ldots,\al_l\})$ is a base of $\pi(R)$.
\par
%{\it{Let $\Pi\cup\{a\}$ be a fundamental-set of $R$, where
%$\Pi=\{\al_0,\al_1\,\ldots,\al_l\}$; and moreover,
%we assume $\al_0$ to be such that}}
%\begin{equation}\label{eqn:eqndefellalzero}
%\begin{array}{l} 
%\mbox{$\{\pi(\al_1)\,\ldots,\pi(\al_l)\}$ is a
%base of $\pi(R)$}  \\
%\mbox{and, if $2\pi_a(\al_0)\in \pi_a(R)$, then $2\pi_a(\al_i)\in \pi_a(R)$ for some $i\ne 0$}
%\end{array}
%\end{equation}
%(cf. Theorem~\ref{theorem:theoremoneone}).
Let $\delta(\Pi)\in\bbbz \Pi$
be such that
\begin{equation}\label{eqn:eqndefelldel}
\delta(\Pi)\in\bbbn \Pi\quad\mbox{and}\quad
\bbbz \delta(\Pi)=(\bbbz \Pi)^0.
\end{equation} 
Then $\pi_a(\delta(\Pi))=\delta(\pi_a(\Pi))$ (see (\ref{eqn:eqndefondelpi})
for $\delta(\pi_a(\Pi))$). \par
Let $\delta=\delta(\Pi)$ be as in 
(\ref{eqn:eqndefelldel}).
By (\ref{eqn:bsofexaffd}), (\ref{eqn:eqnzraff})
and (\ref{eqn:eqndfofbase}),
for $\bbbx\in\{\bbbz,\bbbr\}$, we have
\begin{equation}\label{eqn:eqnzbaseofel}
\left\{\begin{array}{l}
\bbbx R=\bigoplus_{\lambda\in\Pi\cup\{a\}}\bbbx \lambda
=(\bigoplus_{\al\in\Pi\setminus\{\al_0\}}\bbbx \al)\bigoplus\bbbx 
\delta\bigoplus\bbbx a, \\
(\bbbx R)^0 =\bbbx \delta\oplus\bbbx a, \\
R\subset (\bbbxp \Pi\oplus\bbbx a)\cup (\bbbxm \Pi\oplus\bbbx a).
\end{array}\right.
\end{equation} 
%\par Set
%\begin{equation}\label{eqn:defofapiapi}
%A(\Pi):=A(\pi_a(\Pi))
%\end{equation} (see (\ref{eqn;eqndefofapi}) for $A(\pi_a(\Pi))$).

\subsection{Maps $k$ and $g$ }

\begin{lmm} \label{lemma:alS} {\rm{(1)}} For any $\al\in R$, we have
\begin{equation}\label{eqn:eqnalfiber}
(\al+(\bbbz\setminus\{0\})a)\cap R\ne\emptyset.
\end{equation}
\par
{\rm{(2)}}
Let $S$ be a non-empty proper connected subset of $\Pi$. 
Let $\v^S:=\bbbr S\oplus\bbbr a$ and $R^S:=R\cap \v^S$. Then $(R^S,\v^S)$
is a reduced affine root system
(we have assumed $R$ is reduced),
and $(\pi_a(R^S), \v /\bbbr a)$ is an irreducible finite root system
with the base $\pi_a(S)$. In particular,
$\bbbz R^S=\bbbz S\oplus \bbbz k_S a$ for some $k_S\in\bbbn$.

\end{lmm}
\par {\it{Proof.}} {\rm{(1)}}
By (\ref{eqn:eqnzbaseofel}),
$R$ cannot be included in $\bbbz\Pi$.
Hence there exist $\mu\in R$ and
$m\in\bbbz\setminus\{0\}$
such that $\mu\in ma +\bbbz\Pi$. Since $\pi_a(R)$ is an affine root system
and $\pi_a(\Pi)$ is a base of $\pi_a(R)$, by the first equality of (\ref{eqn:affnonrt}),
there exist $\gamma\in\Pi$,
$c\in\{1,2\}$ and $w\in W_\Pi$
such that $w(\mu)=c\gamma+ma$. Notice that
\begin{equation}\label{eqn:eqnaltcma}
R\ni s_\gamma s_{c\gamma+ma}(\gamma)
=s_\gamma(\gamma-(c^{-1}2)(c\gamma+ma))=\gamma-2c^{-1}ma.
\end{equation} (Hence (\ref{eqn:eqnalfiber})
holds for this special $\gamma$.)
Let $\lambda=\gamma-2c^{-1}ma$.
For $\beta\in R$, we have
\begin{equation}\label{eqn:eqmrnisaslambda}
R\ni s_\gamma
s_\lambda(\beta)=s_\gamma(\beta-(\gamma^\vee,\beta)\lambda)
=\beta+(\gamma^\vee,\beta)\cdot 2c^{-1}ma.
\end{equation}
By ({\rm{AX}}5) and (\ref{eqn:eqnzbaseofel}),
by repetition of equations similar to (\ref{eqn:eqmrnisaslambda}),
we see that (\ref{eqn:eqnalfiber}) holds for any
$\al\in R$. \par
{\rm{(2)}} This follows from (1) and (\ref{eqn:eqnzbaseofel}).
\hfill $\Box$
\newline\par
By Lemma~\ref{lemma:alS}~(2), for each $\al\in\Pi$, $R^{\{\al\}}$ is a rank-one 
reduced affine
root system and $\{\pi_a(\al)\}$ is a base of a rank-one irreducible finite
root system $\pi_a(R^{\{\al\}})$.
By Theorem~\ref{theorem:theoremoneone}, we can define maps
\begin{equation}\label{eqn:eqnkgdis}
k:\Pi\rightarrow\bbbn \,\,\mbox{and}\,\,
g:\Pi\rightarrow\{\emptyset,2\bbbz +1 \}
\end{equation}
by
\begin{equation}\label{eqn:eqndefkgpre}
R\cap(\bbbr \al\oplus\bbbr a)=
\bigcup_{\varepsilon\in \{1,-1\}}
((\varepsilon\al+\bbbz k(\al)a)\cup(2\varepsilon\al+g(\al)k(\al)a))
\end{equation} $(\al\in\Pi)$ (
see also (\ref{eqn:eqnzbaseofel})).

Since $\pi_a(R)\setminus 2\pi_a(R)= W_{\pi_a(\Pi)}\cdot \pi_a(\Pi)$
(see Theorem~\ref{theorem:theoremoneone}),
we have
\begin{equation}\label{eqn:eqnexppre}
R=\bigcup_{w\in W_\Pi}(\bigcup_{\al\in\Pi}
((w(\al)+\bbbz k(\al)a)\cup(w(2\al)+g(\al)k(\al)a))).
\end{equation} 
Since
$R$ is determined by $\Pi$, $k$ and $g$, 
\begin{equation}\label{eqn:eqnpreexp}
\mbox{we also denote $R$ by $R(\Pi,k,g)$.}
\end{equation} 

Let $\al\in\Pi$. Let $\al^*:=-\al_0(R^{\{\al\}},\{\al\},-k(\al)a)$.
Then $\al^*=c(\al)\al+k(\al)a$, where
\begin{equation}\label{eqn:eqnca}
c(\al)=\left\{\begin{array}{ll}
1 & \quad\mbox{if $g(\al)=\emptyset$,} \\
2 & \quad\mbox{if $g(\al)=2\bbbz +1$.}
\end{array}\right.
\end{equation} 

Let $\bp:=\{\al,\al^*|\al\in\Pi\}$. Then $|\bp|=2|\Pi|=2(l+1)$. By Thereom~\ref{theorem:theoremoneone}, we have
\begin{equation}\label{eqn:wbpeqw}
R=W_\bp\cdot \bp\,\,\mbox{and}\,\,W=W_\bp
\end{equation} (We have assumed that $R$ is reduced).

Assume $l\geq 2$ (see (\ref{eqn:ranknullitydef})).
Let $\al$, $\beta\in\Pi$ be such that $(\beta^\vee,\al)=-1$. 
Let $\gamma=\al_0(R^{\{\al,\beta \}},\{\al,\beta\},-k(\al)a)$. By Lemma~\ref{lemma:alS}~(2)
and Theorem~\ref{theorem:theoremoneone}, we have $g(\beta)=\emptyset$, $k_{\{\al,\beta\}}=k(\al)$ and see that 
$((\beta^\vee,\al),k(\beta)/k(\al),g(\al))$ for 
the rank-two reduced affine root system $R^{\{\al,\beta\}}$ with a base $\{\al,\beta,\gamma\}$ 
is one of the following.
\begin{equation}\label{eqn:listRab}
\left\{\begin{array}{ll}
(-1,1,\emptyset ) & \mbox{so $R^{\{\al,\beta\}}$ is $A_2^{(1)}$, and $\gamma=-s_\al (\beta^*)$,} \\
(-2,1,\emptyset )& \mbox{so $R^{\{\al,\beta\}}$ is $B_2^{(1)}$, and $\gamma=-s_\al (\beta^*)$,} \\
(-3,1, \emptyset ) & \mbox{so $R^{\{\al,\beta\}}$ is $G_2^{(1)}$, and $\gamma=-s_\beta s_\al (\beta^*)$,} \\
(-2,2, \emptyset ) & \mbox{so $R^{\{\al,\beta\}}$ is $D_3^{(2)}$, and $\gamma=-s_\beta (\al^*)$,} \\
(-3,3, \emptyset) & \mbox{so $R^{\{\al,\beta\}}$ is $D_4^{(3)}$, and $\gamma=-s_\al s_\beta (\al^*)$,} \\
(-2,1, 2\bbbz +1) & \mbox{so $R^{\{\al,\beta\}}$ is $A_4^{(2)}$, and $\gamma=-s_\beta (\al^*)$.}
\end{array}\right.
\end{equation}

\subsection{List of $(\Pi,k,g)$}

\begin{thm} Let $R=R(\Pi,k,g)$ be as in (\ref{eqn:eqnpreexp}).

{\rm{(1)}} Assume $l=1$. Let $\{\al_1,\al_0\}=\Pi$ and assume that $\{\pi(\al_1)\}$
is a base of $\pi(R)$ and that $k(\al_1)\leq k(\al_0)$ if $\{\pi(\al_0)\}$
is also a base of $\pi(R)$.
Then $k(\al_1)=1$ and $((\al_0^\vee,\al_1),k(\al_0),g(\al_0),g(\al_1))$
is exactly one of the followings:
\begin{equation}\begin{array}{l}
(-2,1,\emptyset,\emptyset), \\
(-2,1,\emptyset,2\bbbz +1),\,(-2,1,2\bbbz +1,\emptyset),\,(-2,1,2\bbbz +1,2\bbbz +1), \\
(-2,2,\emptyset,\emptyset),\,(-2,2,2\bbbz +1,\emptyset),\\
(-1,1,\emptyset,\emptyset),\,(-1,1,\emptyset,2\bbbz +1),\,(-1,2,\emptyset,\emptyset),\,(-1,2,\emptyset,2\bbbz +1), \\
(-1,4,\emptyset,\emptyset).
\end{array}
\end{equation}

{\rm{(2)}} Assume $l\geq 2$. Then there exists $R(\Pi,k,g)$ such that $(W_\Pi\cdot \Pi,\bbbr\Pi)$ is a rank-$l$ reduced
affine root system of any type with a base $\Pi$
and $k:\Pi\to\bbbn$ and $g:\Pi\to\{\emptyset,2\bbbz +1\}$ are any maps satisfying the condition that $1\in k(\Pi)$
and
$((\al^\vee,\beta),k(\beta)/k(\al),g(\al))$ is the same as one of {\rm{(\ref{eqn:listRab})}}
for any $\al$, $\beta\in\Pi$ with $(\beta^\vee,\al)=-1$. 
\end{thm}

The statements of this theorem is well-known and, however,
some of $R(\Pi,k,g)$'s are isomorphic
(see \cite[(6.6)]{saito1} 
and \cite[Lists~4.6, 4.25, 4.67, 4.78]{Allisonetal}).
For the case $\l\geq 2$, which of them are isomorphic
can be read off from the statement of Theorem~\ref{theorem:dimnull}. 

\section{Elliptic Lie algebras with rank $\geq 2$}\label{section:ellipticLie}

In this section we assume $R$ is a reduced elliptic root system with rank $\geq 2$, that is, 
$R\cap 2R=\emptyset$, $n=2$ and 
$l\geq 2$
{\rm{(}}see  {\rm{(\ref{eqn:ranknullitydef}))}}. We have assumed the 
rank $l\geq 2$
mainly because
we use the fact 
(\ref{eqn:rktwoafli}) below.
We fix a fundamental-set $\Pi\cup\{a\}$ of $R$.

\subsection{Useful lemma}

The following lemma is useful.

\begin{lmm} \label{lemma:cpcal} Let $\v^\prime$ be a 2-dimensional $\bbbc$-linear space having
a non-degenerate symmetric 
bilinear form $(\,,\,):\v^\prime\times\v^\prime\to\bbbc$.
Let $\gamma_1$, $\gamma_2\in (\v^\prime)^\times$.
Let $\fa$ be a Lie algebra over $\bbbc$ generated by ${\bar h}_\gamma$ ($\gamma\in\v^\prime$),
${\bar E}_1$, ${\bar E}_2$, ${\bar F}_1$, ${\bar F}_2$ and satisfying the equations 
${\bar h}_{x\gamma+x^\prime\gamma^\prime}=x{\bar h}_\gamma+x^\prime{\bar h}_{\gamma^\prime}$, 
$[{\bar h}_\gamma,{\bar h}_{\gamma^\prime}]=0$,
$[{\bar h}_\gamma,{\bar E}_i]=(\gamma,\gamma_i){\bar E}_i$,
$[{\bar h}_\gamma,{\bar F}_i]=-(\gamma,\gamma_i){\bar F}_i$,
and $[{\bar E}_i,{\bar F}_i]=\delta_{ij} {\bar h}_{\gamma_i^\vee}$,
for $x$, $x^\prime\in\bbbc$, $\gamma$, $\gamma^\prime\in\v^\prime$,
and $i\in\{1,2\}$. 

{\rm{(1)}} For $k\in\bbbn$, we have
\begin{equation}\label{eqn:adfone}
\begin{array}{l}
[\ad ({\bar E}_1)^k({\bar E}_2),\ad ({\bar F}_1)^k({\bar F}_2)] \\
\quad =
k!(\prod_{m=1}^{k-1}((\gamma_1^\vee,\gamma_2)+m))(k(\gamma_1,\gamma_2^\vee){\bar h}_{\gamma_1^\vee}+(\gamma_1^\vee,\gamma_2){\bar h}_{\gamma_2^\vee}).
\end{array}
\end{equation} 

{\rm{(2)}} Let $m:=(\gamma_1^\vee,\gamma_2)$. Assume $m\in\bbbzm$.
Assume that   
${\bar h}_{\gamma_1^\vee}$ and ${\bar h}_{\gamma_2^\vee}$
are linearly independent.
Assume $\ad ({\bar E}_1)^r({\bar E}_2)=\ad ({\bar F}_1)^r({\bar F}_2)=0$ for some $r\in\bbbn$. Let
\begin{equation}\label{eqn:adop}
{\bar n}=n({\bar E}_1,{\bar F}_1):=\exp(\ad {\bar E}_1)\exp(-\ad {\bar F}_1)\exp(\ad {\bar E}_1).
\end{equation} 
Then we have 
\begin{equation}\label{eqn:adftwo}\begin{array}{l}
\ad ({\bar E}_1)^{1-m}({\bar E}_2)=\ad ({\bar F}_1)^{1-m}({\bar F}_2)=0, \\
{\bar n}({\bar h}_\gamma )={\bar h}_\gamma-(\gamma_1,\gamma){\bar h}_{\gamma_1^\vee},\,
{\bar n}({\bar E}_1)=-{\bar E}_1,\,{\bar n}({\bar F}_1)=-{\bar F}_1, \\
{\bar n}((\ad {\bar E}_1)^i{\bar E}_2)={\frac {(-1)^ii!} {(-m-i)!}}(\ad {\bar E}_1)^{-m-i}{\bar E}_2\ne 0, \\
{\bar n}((\ad {\bar F}_1)^i{\bar F}_2)={\frac {(-1)^{m-i}i!} {(-m-i)!}}(\ad {\bar F}_1)^{-m-i}{\bar F}_2\ne 0,
\end{array}
\end{equation} for $0\leq i\leq -m$ and $\gamma\in\v^\prime$.
\end{lmm}

We can get (\ref{eqn:adfone}) directly and get (\ref{eqn:adftwo}) by using a representation theory
of ${\rm{sl}}_2$.

\subsection{Definition of elliptic Lie algebras with rank $\geq 2$}

Let 
$\a:=\{(\al,\beta)\in \Pi\times\Pi\,
|\,(\al,\beta^\vee)=-1\}$.
Let $\b:=\bp\cup(-\bp)$, and
$\btp:=\{(\mu,\nu)\in\b\times\b\,|\,\mu\ne\nu\ne -\mu\}$.
For $(\mu,\nu)\in\btp$, let $x_{\mu,\nu}=1-((\mu^\vee,\nu)-|(\mu^\vee,\nu)|)/2$.
Let $\vc=\bbbc\otimes_\bbbr\v$, so $\vc$ is a $l+2$-dimensional $\bbbc$-linear space. 
We identify $\v$ with the $\bbbr$-linear subspace $1\otimes \v$ of $\vc$.
We say that a map $\omega:\a\to\bbbct$ is a {\it{tuning}} if $\omega (\al,\beta)\omega (\beta,\al)=1$
whenever $(\al^\vee,\beta)=-1$. Denote $\omega_1$ by the tuning with $\omega_1 (\al,\beta)=1$ for
all $(\al,\beta)\in\a$, and moreover, if $W_\Pi\cdot \Pi$ is $A_l^{(1)}$, then for $q\in\bbbct$, 
denote $\omega_q$ by the tuning with $\omega (\al_i,\al_{i+1})=1$ ($0\leq i\leq l$) and $\omega (\al_l,\al_0)=q$,
where the numbering of the elements of $\Pi$ is the same as that of the Dynkin diagram
of $A_l^{(1)}$ in Subsection~\ref{subsection:superaffdynkin}.

\begin{dfn}
\label{definition:defofellsuper}
Let
$\fg ^\omega=\fg (\Pi,k,g,\omega)$ be the
Lie algebra over $\bbbc$ defined by generators:
\begin{equation}\nonumber
h_\sigma\,\,(\sigma\in\vc), \quad
E_\mu\,\,(\mu\in\b),
\end{equation} and relations: \newline\par
({\rm {SR}}1) \quad
$xh_\sigma+yh_\tau=h_{x\sigma+y\tau}$
if $x,\,y\in\bbbc$ and $\sigma,\,\tau\in\vc$, \par
({\rm {SR}}2) \quad $[h_\sigma,h_\tau]=0$ if
$\sigma,\,\tau\in\vc$, \par
({\rm {SR}}3) \quad
$[h_\sigma,E_\mu]=(\sigma,\mu)E_\mu$
if $\sigma\in\vc$ and $\mu\in\b$, \par
({\rm {SR}}4) \quad
$[E_\mu,E_{-\mu}]=h_{\mu^\vee}$
if $\mu\in\bp$, \par
({\rm {SR}}5) \quad
$(\ad E_\mu)^{x_{\mu,\nu}}E_\nu=0$
if $(\mu,\nu)\in \btp$,  \par
({\rm {SR}}6) \quad
$c(\al)(\ad E_{\al^*})^\kbka E_\beta=\omega(\al,\beta)({\rm
ad}E_\al)^{c(\al)\kbka }E_{\beta^*}$ 
if $(\al,\beta)\in{\cal A}$, \par
({\rm {SR}}7) \quad
$(-1)^{c(\al)+1}c(\al)({\rm
ad}E_{-\al^*})^\kbka E_{-\beta}= {\frac 1 {\omega(\al,\beta)}}({\rm
ad}E_{-\al})^{c(\al)\kbka }E_{-\beta^*}$
if $(\al,\beta)\in{\cal A}$, \par
({\rm {SR}}8) \quad
$(\ad E_\al)^i(\ad E_{\al^*})^{\kbka -i}E_\beta=0$
if $(\al,\beta)\in{\cal A}$ and $1\leq i\leq \kbka -1$, \par
({\rm {SR}}9) \quad
$(\ad E_{-\al})^i({\rm
ad}E_{-\al^*})^{\kbka -i}E_{-\beta}=0$
if $(\al,\beta)\in{\cal A}$ and $1\leq i\leq \kbka -1$.
\newline\par
We call $\fg (\Pi,k,g,\omega)$ an {\it{elliptic Lie algebra}},
see Introduction. Let
$\fg =\fg (\Pi,k,g):=\fg ^{\omega_1}$. 
\end{dfn}

We have
\begin{lmm}\label{lemma:tuning}
If $W_\Pi\cdot \Pi$ is not $A_l^{(1)}$ (resp. is $A_l^{(1)}$),
then there is an isomorphism $\varphi$ from $\fg^\omega$ to $\fg$ (resp. to $\fg^{\omega_q}$ for some $q\in\bbbct$)
such that $\varphi(h_\sigma)=h_\sigma$ $(\sigma\in\vc)$ and $\varphi(E_\mu)\in\bbbct E_\mu$
$(\mu\in\b)$.
\end{lmm}

{\it{Proof.}} Using (\ref{eqn:adfone}), we can modify ({\rm {SR}}6-7) by 
taking non-zero scalar products of $E_\mu$'s. \hfill
$\Box$
\newline\par
Let $\fh^\omega=\fh^\omega(\Pi,k,g,\omega):=\{h_\sigma\in\fg ^\omega|\sigma\in\vc\}$,
and $\fh=\fh (\Pi,k,g):=\fh^{\omega_1}$.

Since all equations in ({\rm {SR}}1-9) are $\bbbz R$-homogeneous, where $R=R(\Pi,k,g)$, we can regard
$\fg^\omega$ as
the $\bbbz R$-graded Lie algebra $\fg^\omega=\oplus_{\sigma\in \bbbz R}\fg^\omega_\sigma$
(that is $[\fg^\omega_\sigma,\fg^\omega_{\sigma^\prime}]\subset \fg^\omega_{\sigma+\sigma^\prime}$)
such that $E_\mu\in \fg^\omega_\mu$ for all $\mu\in\b$.
Note $\fh^\omega\subset \fg^\omega_0$. For each $\mu\in\bp$, we can define $n_\mu$ to be $n(E_\mu,E_{-\mu})$ 
(see (\ref{eqn:adop})) as 
an automorphism of $\fg^\omega$, so $n_\mu(\fg^\omega_\sigma)= \fg^\omega_{s_\mu(\sigma)}$.
Let $\r^\omega=\{\sigma\in \bbbz R|\dim \fg^\omega_\sigma\ne 0\}$. Then we have
\begin{equation}\label{eqn:wro}
W_\bp\cdot \r^\omega=\r^\omega.
\end{equation}
Let $S$ a non-empty proper connected subset of $\Pi$. Let $\fg^{\omega,S}$
be the Lie algebra over $\bbbc$ defined by the generators
$h_\sigma$ ($\sigma\in\bbbc S\oplus\bbbc a$), $E_{\pm\al}$, $E_{\pm\al^*}$ ($\al\in S$)
and the same relations as those in ({\rm{SR}}1-9). Let $\iota^{\omega,S}:\fg^{\omega,S}\to
\fg^\omega$ be the homomorphism sending the generators to those denoted by the same symbols. 
Let $\fg^{\omega,S}_\sigma=(\iota^{\omega,S})^{-1}(\fg^\omega_\sigma)$
for $\sigma\in\bbbz R^S$, so  $\fg^{\omega,S}=\oplus_{\sigma\in \bbbz R^S}\fg^{\omega,S}_\sigma$.
Let $\fg^S=\fg^{\omega_1,S}$, and $\fg^S_\sigma=\fg^{\omega_1,S}_\sigma$.
Let $\r^{\omega,S}=\{\sigma\in \bbbz R^S|\dim \fg^{\omega,S}_\sigma\ne 0\}$. 

Let $\al\in\Pi$. Then $\fg^{\omega,\{\al\}}=\fg^{\{\al\}}$,
since $\fg^{\omega,\{\al\}}$ is defined by using ({\rm{SR}}1-5).
By Serre's relations
({\rm{SR}}1-5), 
 $\fg^{\omega,\{\al\}}$ is (the derived  algebra of) an
affine Lie algebra with $\r^{\omega,\{\al\}}=R^{\{\al\}}\cup\bbbz k(\al)a$,
where the affine root system $R^{\{\al\}}$
is $A_1^{(1)}$ or $A_2^{(1)}$. Hence 
$\dim \fg^{\omega,\{\al\}}_0=2$, and 
$\dim \fg^{\omega,\{\al\}}_\lambda=1$ 
($\lambda\in\r^{\omega,\{\al\}}\setminus\{0\}$). Note 
$\r^{\omega,\{\al\}}\setminus\{0\}=R^{\{\al\}}\cup\bbbzt k(\al)a$.

\begin{lmm} \label{lemma:chiainj} There is a homomorphism $\chi^\omega$ from $\fg^\omega$ to a Lie algebra $\fb^\omega$ 
such that $\dim\chi^\omega(\fh^\omega)=l+2$, 
$\dim\chi^\omega(\iota^{\omega,\{\al\}}(\fg^{\omega,\{\al\}}_\lambda))=1$
for all $\al\in\Pi$ and all $\lambda\in R^{\{\al\}}\cup\bbbzt k(\al)a$, and 
\begin{equation}\label{eqn:chipr} 
\begin{array}{l}
\chi^\omega(\fh^\omega+\sum_{\al\in\Pi}\sum_{\lambda\in R^{\{\al\}}\cup\bbbzt k(\al)a}
\iota^{\omega,\{\al\}}(\fg^{\omega,\{\al\}}_\lambda)) \\
\quad =\chi^\omega(\fh^\omega)\oplus
\bigoplus_{\al\in\Pi}\bigoplus_{\lambda\in R^{\{\al\}}\cup\bbbzt k(\al)a}
\chi^\omega(\iota^{\omega,\{\al\}}(\fg^{\omega,\{\al\}}_\lambda)).
\end{array} 
\end{equation}
\end{lmm}

(If $\omega=\omega_1$, then $\fb^\omega$ is given as an `affinization' 
$\mathfrak{a}\otimes\bbbc [t,t^{-1}]\oplus\bbbc c$
of (the derived algebra of) an affine Lie algebra $\mathfrak{a}$, see \cite[Proposition~3.1]{yamane2}.)

{\it{Proof.}} If $\omega=\omega_1$, then we can define $\chi=\chi^{\omega_1}$ in a
way entirely similar to that of \cite[Proposition~3.1]{yamane2}, inspired by
so-called 
an `unfolding process' of a Dynkin diagram of a reduced affine root system,
and we see by checking each case directly that such $\chi$ has the property (\ref{eqn:chipr}). 
The existence of a $\chi^{\omega_q}$ is well-known (see \cite{BermanGaoKrylyuk}).
Then this lemma follows from Lemma~\ref{lemma:tuning}.
\hfill $\Box$ 

For each $\al\in\Pi$, let $[R^{\{\al\}}]^+:=R^{\{\al\}}\cap(\bbbn\al + \bbbz k(\al)a)$,
and $[R^{\{\al\}}]^-:=-[R^{\{\al\}}]^+$.
Note that $R^{\{\al\}}=[R^{\{\al\}}]^+\cup [R^{\{\al\}}]^-$.

\begin{lmm} \label{lemma:chiainjd} 
For each $(\al,\beta)\in\a$, 
\begin{equation}\label{eqn:rktwoafli}
\begin{array}{l}
\mbox{$\fg^{\omega,\{\al,\beta\}}$ is (the derived algebra of) an affine Lie algebra} \\
\mbox{with the affine root system $R^{\{\al,\beta\}}$,}
\end{array}
\end{equation}
which implies $\r^{\omega,\{\al,\beta\}}=R^{\{\al,\beta\}}\cup\bbbz k(\al)a$.
In particular, for each $(\al^\prime,\beta^\prime)\in\Pi\times\Pi$ with 
$\al^\prime\ne \beta^\prime$,  we have
\begin{equation}\label{eqn:pmapmbwt}
[\iota^{\omega,\{\al^\prime\}}(\fg^{\omega,\{\al^\prime\}}_\lambda), 
\iota^{\omega,\{\beta^\prime\}}(\fg^{\omega,\{\beta^\prime\}}_\mu)]=0
\end{equation} for all $(\lambda,\mu)\in ([R^{\{\al^\prime\}}]^+ \times [R^{\{\beta^\prime\}}]^-)\cup
([R^{\{\al^\prime\}}]^- \times [R^{\{\beta^\prime\}}]^+)$.
\end{lmm}

{\it{Proof.}} Note first that $h_\al$, $h_\beta$ and $h_a$ are linearly independent in 
$\fg^{\omega,\{\al,\beta\}}$, which follows from
Lemma~\ref{lemma:chiainj}.
Let $\gamma\in R^{\{\al,\beta \}}$ be as in (\ref{eqn:listRab}). 
If $\gamma$ is expressed as 
$-s_{\gamma_1}\ldots s_{\gamma_{r-1}}(\gamma_r^*)$ in 
(\ref{eqn:listRab}) with $\gamma_i\in\{\al,\beta\}$, then we let $E_{\pm\gamma}:=n_{\gamma_1}\ldots n_{\gamma_{r-1}}(E_{\mp \gamma_r^*})
\in\fg^{\omega,\{\al,\beta \}}_{\pm\gamma}$. Let $\gamma_{r+1}\in\{\al,\beta\}\setminus\{\gamma_r\}$.
By ({\rm{SR}}6-7) and (\ref{eqn:adftwo}), 
we have $n_{\pm \gamma_r^*}(E_{\pm \gamma_{r+1}})=n_{\pm \gamma_r}(E_{\pm \gamma_{r+1}^*})$.
Hence $\fg^{\omega,\{\al,\beta \}}$ is generated by $E_{\pm \al}$, $E_{\pm \beta}$ and 
$E_{\pm \gamma}$. We show 
\begin{equation}\label{eqn:pmapmbmpg}
[E_{\pm \al},E_{\mp \gamma}]=[E_{\pm \beta},E_{\mp \gamma}]=0.
\end{equation} If $R^{\{\al,\beta\}}\ne A_4^{(2)}$, we have this in the same way as in 
\cite[\S 2.3]{yamane2}.
Assume $R^{\{\al,\beta\}}= A_4^{(2)}$. We write $X \sim Y$ if $X\in\bbbct Y$.
By (\ref{eqn:adftwo}) and ({\rm{SR}}6), 
\begin{equation}
E_{-\gamma}\sim  [E_\beta,[E_\beta, E_{\al^*}]]
\sim  [E_\beta,[E_\al, [E_\al, E_{\beta^*}]]]
\label{eqn:twosubcasesix}
\end{equation}
Then $[E_\beta,E_{-\gamma}]=0$
follows from ({\rm{SR}}5).
We have
\begin{eqnarray*}
                  [E_{-\gamma},E_\al ] &\sim &
[[E_\beta,[E_\al, [E_\al, E_{\beta^*}]]],E_\al ]
\quad\mbox{(by (\ref{eqn:twosubcasesix}))} \\
&\sim & [[E_\beta, E_\al ],[E_\al, [E_\al, E_{\beta^*}]]]]
\quad\mbox{(by ({\rm{SR}}5))} \\
&\sim & [[E_\beta, E_\al ],[E_\beta, E_{\al^*}]]]
\quad\mbox{(by ({\rm{SR}}6))} \\
&\sim & n_\beta([E_\al ,[E_\beta, E_{\al^*}]])\quad\mbox{(by (\ref{eqn:adftwo}))} \\
&\sim & n_\beta([E_\al ,[E_\al, [E_\al, E_{\beta^*}]]])
\quad\mbox{(by ({\rm{SR}}6))} \\
&= & 0\quad\mbox{(by ({\rm{SR}}5))}.
\end{eqnarray*} The remaining equalities of (\ref{eqn:pmapmbmpg})
can be shown similarly.
Hence by (\ref{eqn:adftwo}) and ({\rm{SR}}5), the above generators satisfy Serre's relations. 
Hence (\ref{eqn:rktwoafli}) holds, as desired. \hfill $\Box$
\newline\par
For $i\in\bbbn$, let $(\fn^{\omega,\pm})^{(i)}$ be the $\bbbc$-linear subspaces 
of $\fg^\omega$ defined by
$(\fn^{\omega,\pm})^{(1)}:=\oplus_{\al\in\Pi}\oplus_{\lambda\in[R^{\{\al\}}]^\pm}\iota^{\omega,\{\al\}}(\fg^{\omega,\{\al\}}_\lambda)$
(see Lemme~\ref{lemma:chiainj}),
and $(\fn^{\omega,\pm})^{(i)}:=[(\fn^{\omega,\pm})^{(1)},(\fn^{\omega,\pm})^{(i-1)}]$
inductively for $i\geq 2$.
Let $\fn^{\omega,\pm}$ be the two Lie subalgebras of $\fg^\omega$ defined by 
$\fn^{\omega,\pm}
:=\sum_{i=1}^\infty (\fn^{\omega,\pm})^{(i)}$.
Let $\fn^{\omega,\pm}_\sigma=\fg^\omega_\sigma\cap \fn^{\omega,\pm}$. 
Then
$\fn^{\omega,\pm}=\oplus_{\sigma\in (\bbbz_{\pm} \Pi\oplus\bbbz a)\setminus\bbbz a}\fn^{\omega,\pm}_\sigma$.
For each $\al\in\Pi$,
since $\iota^{\omega,\{\al\}}$ is a Lie algebra homomorphism
(preserving $\bbbz\Pi\oplus\bbbz a$-grading),
we have $\fn^{\omega,\pm}_\mu=
\fn^{\omega,\pm}_\mu\cap(\fn^{\omega,\pm})^{(1)}
=\iota^{\omega,\{\al\}}(\fg^{\omega,\{\al\}}_\mu)$
for all $\mu\in(\bbbz_\pm \al\oplus\bbbz a)\setminus \bbbz a$.
Moreover, by (\ref{eqn:pmapmbwt}), we have 
\begin{equation}[(\fn^{\omega,+})^{(1)},(\fn^{\omega,-})^{(1)}]\subset
(\fn^{\omega,+})^{(1)} + (\fn^{\omega,-})^{(1)} +
\sum_{\al\in\Pi}
\sum_{\sigma\in \bbbz k(\al)a}
\iota^{\omega,\{\al\}}(\fg^{\omega,\{\al\}}_\sigma ).
\end{equation} 
Hence by Lemma~\ref{lemma:chiainj} and (\ref{eqn:rktwoafli}), we have 
\begin{equation}\label{eqn:elltri}
\fg^\omega = \fh^\omega\oplus \fn^{\omega,+}\oplus \fn^{\omega,-}\oplus(\bigoplus_{\al\in\Pi}
\bigoplus_{\sigma\in \bbbzt k(\al)a}
\iota^{\omega,\{\al\}}(\fg^{\omega,\{\al\}}_\sigma )),
\end{equation} $\dim\fh^\omega=l+2$, and $\dim \fn^{\omega,\pm}_\lambda=\dim \iota^{\omega,\{\al\}}(\fg^{\omega,\{\al\}}_\sigma )=1$
for $\al\in\Pi$, $\lambda\in[R^{\{\al\}}]^\pm$ and $\sigma\in \bbbzt k(\al)a$.
By (\ref{eqn:affnonrt}), we have
\begin{equation}\label{eqn:ellnonrt}
\left\{\begin{array}{l}
R=W_\Pi\cdot \bigcup_{\al\in\Pi}[R^{\{\al\}}]^+,  \\
(\bbbz R)^\times\setminus R \\
=W_\Pi\cdot (\bigcup_{\al\in\Pi}(\bbbn\al \oplus\bbbz a)\setminus [R^{\{\al\}}]^+)\cup((\bbbz R)^\times\setminus(\bbbzp\Pi\cup\bbbzm\Pi)\oplus\bbbz a).
\end{array}\right.
\end{equation}
Then by (\ref{eqn:wro}), using a standard argument as in \cite{kac1}, \cite{saitoyoshii},
together with the automorphisms $n_\mu$ ($\mu\in\bp$), we have

\begin{thm}\label{theorem:mmain} We have 
$(\r^\omega)^\times = R$, $\dim \fg^\omega_\mu=1$ $(\mu\in R)$,
$\fg^\omega_0=\fh^\omega$, $\dim\fh^\omega=l+2$, 
$(\r^\omega)^0 \subset \bbbz \delta\oplus \bbbz a$,
and $\dim \fg^\omega_{ma}=|\{\al\in\Pi|m\in \bbbz k(\al)\}|$
 $(m\in\bbbzt)$. 
\end{thm}

By the following theorem, we can compute $\dim \fg^\omega_\lambda$
for $\lambda\in \bbbz\delta \oplus \bbbz a$. 

\begin{thm}\label{theorem:mainuni} 
Let $\Pi^\prime\cup\{a^\prime\}$ be a fundamental-set of $R$.
Then there exist a tuning $\eta$ for $\Pi^\prime\cup\{a^\prime\}$
and an isomorphism $f:\fg(\Pi^\prime,k^\prime,g^\prime,\eta)
\to \fg^\omega$ such that $f(\fg^{\prime,\eta}_\lambda)=\fg^\omega_\lambda$
for all $\lambda\in \bbbz\Pi\oplus\bbbz a$, where $\fg^{\prime,\eta}:=\fg(\Pi^\prime,k^\prime,g^\prime,\eta)$.
In particular, we have
\begin{equation}\label{eqn:unieq} 
\dim \fg^\omega_{ma^\prime}=|\{\al^\prime\in\Pi^\prime|m\in \bbbz k^\prime(\al^\prime)\}|
\,\,\mbox{for $m\in\bbbzt$}.
\end{equation}
\end{thm}

{\it{Proof.}} Let $\bp^\prime=\{\al^\prime,(\al^\prime)^*|\al\in\Pi^\prime\}$
and $\b^\prime=\bp^\prime\cup -\bp^\prime$.
By ({\rm{SR}}1-9), Theorem~\ref{theorem:mmain} and (\ref{eqn:adftwo}),
for some $\eta$, we have a homomorphism $f$ of the statement such that
$f(\fg^{\prime,\eta}_{\mu^\prime})=\fg^\omega_{\mu^\prime}$ for all 
$\mu^\prime\in \b^\prime$.
Since $\fg^{\prime,\eta}$ is generated by $\fg^{\prime,\eta}_{\mu^\prime}$ ($\mu^\prime\in \b^\prime$), 
we have  $f(\fg^{\prime,\eta}_\lambda)\subset \fg^\omega_\lambda$
for all $\lambda\in\bbbz R=\bbbz\Pi^\prime\oplus \bbbz a^\prime$.
Since $R=W_{\bp^\prime}\cdot \bp^\prime$ by (\ref{eqn:wbpeqw}), using $n(E_{\mu^\prime},E_{-\mu^\prime})
\in{\rm{Aut}}(\fg^{\prime,\eta})$ ($\mu^\prime\in \b^\prime$), by Theorem~\ref{theorem:mmain}, we have 
$f(\fg^{\prime,\eta}_\beta ) =\fg^\omega_\beta$ for all $\beta\in R$. Since
$E_\mu\in f(\fg^{\prime,\eta})$ for all $\mu\in \b$, we have $f(\fg^{\prime,\eta})=\fg^\omega$,
so $f(\fg^{\prime,\eta}_\lambda)= \fg^\omega_\lambda$
for all $\lambda\in\bbbz R$.
By the same argument, for some tuning $\omega^\prime$ for $\Pi\cup \{a\}$, we have an epimorphism
$f^\prime:\fg^{\omega^\prime}=\fg(\Pi,k,g,\omega^\prime)\to \fg^{\prime,\eta}$
such that $f^\prime(\fg^{\omega^\prime}_\lambda)= \fg^{\prime,\eta}_\lambda$
for all $\lambda\in\bbbz R$. Hence $\dim \fg^{\omega^\prime}_\lambda \geq 
\dim \fg^\omega_\lambda$ for all $\bbbz R$, so 
$(\r^\omega)^0 \subset (\r^{\omega^\prime})^0$.
Assume that $W_\Pi\cdot \Pi$ is not $A_l^{(1)}$. By Lemma~\ref{lemma:tuning},
we have $\dim \fg^{\omega^\prime}_\lambda = \dim \fg_\lambda =
\dim \fg^\omega_\lambda$ for all $\lambda\in\bbbz R$, so 
$(\r^\omega)^0 = (\r^{\omega^\prime})^0$.
Hence $f\circ f^\prime$ is an isomorphism, so is $f$. 
Assume that $W_\Pi\cdot \Pi$ is $A_l^{(1)}$.  Assume
$\varphi:\fg(\Pi,k,g,\omega_{q_1})\to \fg(\Pi,k,g,\omega_{q_2})$ is an epimorphism
such that $\varphi(\fg(\Pi,k,g,\omega_{q_1})_\lambda)=\fg(\Pi,k,g,\omega_{q_2})_\lambda$
for all $\lambda\in\bbbz R$. 
For $\gamma\in\bp$, let $c_\gamma\in\bbbct$ be such that $\varphi(E_\gamma)
=c_\gamma E_\gamma$ ($E_\gamma\ne 0$ by Lemma~\ref{lemma:chiainj}). For $\al\in\Pi$, let 
$d_\al=c_\al/c_{\al^*}$.
By ({\rm{SR}6), we have $\omega_{q_2}(\al,\beta)=\omega_{q_1}(\al,\beta)d_\al/d_\beta$
(the element of ({\rm{SR}6) is not zero by Lemma~\ref{lemma:chiainj}
and (\ref{eqn:adfone})).
Hence $d_{\al_i}=d_{\al_{i+1}}$ for $0\leq i\leq l$. Since
$\omega_{q_2}(\al_l,\al_0)=\omega_{q_1}(\al_l,\al_0)$,
we have $q_1=q_2$.
Then by the same argument as above, 
we conclude that 
$f$ is an isomorphism. 

The last statement follows from Theorem~\ref{theorem:mmain}. 
\hfill $\Box$

\section{List of $\dim\fg_{m\delta+ra}$}

In this section we use the notation as follows. For a $\bbbz$-module $X$, $r\in\bbbz$
and $x$, $y\in X$, let
$x \equiv_r y$ means $x-y\in rX$. 
Recall that $l=|\Pi|-1\geq 2$, and
see Subsection~\ref{subsection:superaffdynkin} 
for the numbering of the elements $\al_i$ ($0\leq i\leq l$)
of $\Pi$.
Let $\delta=\delta(\Pi)$.
Fix $\gamma_1\in\PIsh\setminus\{\al_0\}$. 
Fix $\gamma_2\in\PIlg\setminus\{\al_0\}$ 
if $\Rlg\ne\emptyset$. 
Let $M:=\bbbz\delta\oplus \bbbz a$.
We also denote $m\delta +r a\in M$ with $m$, $r\in \bbbz$
by ${m \brack r}$. 
Let $R=R(\Pi,k,g)$ be as in (\ref{eqn:eqnpreexp}). 
Let 
$\Lsh$, $\Llg$ and $\Lex$ be the subsets of $M$ 
such that $\gamma_1+\Lsh=R\cap(\gamma_1+M)$, 
$\gamma_2+\Llg=R\cap(\gamma_2+M)$
(if $\Rlg\ne\emptyset$), and 
$2\gamma_1+\Lex=R\cap(2\gamma_1+M)$ (if $\Rex\ne\emptyset$). 
Let $\Pi^\prime:=\Pi\setminus\{\al_0\}$, so $\pi(\Pi^\prime)$ 
is a base of $\pi(R)$.
By Lemma~\ref{lemma:lemmatheta}, we have
$\Rsh=W_{\Pi^\prime}\cdot \gamma_1+\Lsh$, $\Rlg=W_{\Pi^\prime}\cdot \gamma_2+\Llg$
and $\Rex=W_{\Pi^\prime}\cdot 2\gamma_1+\Lex$.
Let $\fg ^\omega:=\fg (\Pi,k,g,\omega)$,
and $\fg :=\fg ^{\omega_1}$. 

\begin{rem}\label{remark:mod} (Due to Kaiming Zhao) 
Here we would like to mention that a map 
from $M$ to $\{0,1,\ldots,t-1\}$ which is periodic 
modulo $t$ on any line in $M$
is not necessarily meant to be periodic 
modulo $tM$. This indicates that we have to be very careful when 
calculating $\dim\fg^\omega_{m\delta +ra}$ 
because (\ref{eqn:unieq}) does not immediately imply
that $\dim\fg^\omega_{m\delta +ra}$ is periodic, although we finally see that this is true.

Let $f:M\to\bbbzp$ be a map such that $m\bbbz+r\bbbz=f({m \brack r})\bbbz$,
where $f({m \brack r})$ is a ${\rm{g.c.d.}}$ of $m$ and $r$ if
${m \brack r}\ne {0 \brack 0}$.
By definition, $f(h{m \brack r})=h\cdot f({m \brack r})$ for all $h\in\bbbz$
and all ${m \brack r}\in M$.
Let $t\in\bbbn$ be such that $t\geq 2$. 
Define the map $f_t:M\to\{0,1,\ldots,t-1\}$ by 
$f_t({m \brack r})\equiv_tf({m \brack r})$.
Then $f_t((h_1t+h_2){m \brack r})=f_t(h_2{m \brack r})$
for all $h_1\in\bbbz$, all $h_2\in\{0,1,\ldots,t-1\}$
and all ${m \brack r}\in M$.
 Now assume that $t=25$ and ${m \brack r}={40 \brack 200}$.
Then $f({m \brack r})=40$ and $f({m+t \brack r})=5$.
Hence $f_t({m \brack r})=15 \ne 5=f_t({m+t \brack r})$, as desired.

\end{rem}

%\subsection{General examples} 
%By the same argument as in the previous subsection, we can see the following.

Now we have the following theorem.

\begin{thm} \label{theorem:dimnull}
Assume $\fg^\omega=\fg$ if $W_{\Pi}\cdot \Pi$ is not $A_l^{(1)}$ (see Lemma~\ref{lemma:tuning}).
Then $\dim\fg^\omega_\sigma$ with $\sigma\in M\setminus\{0\}$ are listed below.

{\rm{(1)}} Assume that $W_{\Pi}\cdot \Pi$ is $X_l^{(1)}$ with $X=A,\ldots,G$, and $k(\al)=1$ and $g(\al)=\emptyset$
for all $\al\in\Pi$, so 
$\Lsh=M$, $\Rex=\emptyset$, and $\Llg=M$ if $\Rlg\ne\emptyset$ (so $X=B,C,F$ or $G$).
Then we have $\dim\fg ^\omega_\sigma=l+1$ for all 
$\sigma\in M\setminus\{0\}$.

{\rm{(2)}} Assume $W_{\Pi}\cdot \Pi$ is $X_l^{(1)}$ with $X=B,C,F$ or $G$.
Let $r=(\gamma_2,\gamma_2)/(\gamma_1,\gamma_1)$. 
Assume that
$k(\al)=(\al,\al)/(\gamma_1,\gamma_1)$ and $g(\al)=\emptyset$ for all $\al\in\Pi$,
so $\Lsh=M$, $\Llg=\bbbz\delta\oplus \bbbz ra$, and $\Rex=\emptyset$.  
Then we have $\dim\fg_{\sigma_1}=l+1$ for all 
$\sigma_1\in \Llg\setminus\{0\}$, and $\dim\fg_{\sigma_2}=|\PIsh|$ for all 
$\sigma_2\in M\setminus \Llg$.
(This $R$ is isomorphic to $R(\Pi_1,k_1,g_1)$ for which
$W_{\Pi_1}\cdot \Pi_1$ is $D_{l+1}^{(2)}$, $A_{2l-1}^{(2)}$, $E_6^{(2)}$ ($l=4$), or $D_4^{(3)}$ ($l=2$)
respectively, and $k_1(\al )=1$, $g_1(\al )=\emptyset$
($\al\in\Pi$).)

{\rm{(3)}} Assume $W_{\Pi}\cdot \Pi$ is $D_{l+1}^{(2)}$, $A_{2l-1}^{(2)}$, $E_6^{(2)}$ ($l=4$), or $D_4^{(3)}$ ($l=2$).
Let $r=(\gamma_2,\gamma_2)/(\gamma_1,\gamma_1)$. 
Assume that
$k(\al)=(\al,\al)/(\gamma_1,\gamma_1)$ and $g(\al)=\emptyset$ for all $\al\in\Pi$,
so $\Lsh=M$, $\Llg=rM$, and $\Rex=\emptyset$.  
Then we have $\dim\fg_{\sigma_1}=l+1$ for all 
$\sigma_1\in \Llg\setminus\{0\}$, and $\dim\fg_{\sigma_2}=|\PIsh|$ for all 
$\sigma_2\in M\setminus rM$.

{\rm{(4)}} Assume $W_{\Pi}\cdot \Pi$ is $D_{l+1}^{(2)}$, and 
$k(\al_0)=2$, $k(\al_1)=1$, $k(\beta)=2$ ($\beta\in\PIlg$), $g(\al)=\emptyset$ ($\al\in\Pi$),
so $\Lsh=\{0,\delta,a\}+M$, $\Llg=2M$, and $\Rex=\emptyset$.
Then we have $\dim\fg_{\sigma_1}=l+1$ for all 
$\sigma_1\in 2M\setminus\{0\}$, and $\dim\fg_{\sigma_2}=1$ for all 
$\sigma_2\in M\setminus 2M$.

{\rm{(5)}} Assume $W_{\Pi}\cdot \Pi$ is $D_{l+1}^{(2)}$, and 
$k(\al_0)=2$, $g(\al_0)=2\bbbz +1$,  $k(\al_1)=1$,  $g(\al_1 )=\emptyset$, $k(\beta)=2$,
$g(\beta)=\emptyset$  ($\beta\in\PIlg$),
so $\Lsh=\{0,\delta,a\}+M$, $\Llg=2M$ and ${\frac 1 2}\Lex=\delta+a+2M$.
Then we have $\dim\fg_{\sigma_1}=l+1$ for all 
$\sigma_1\in 2M\setminus\{0\}$, and $\dim\fg_{\sigma_2}=1$ for all 
$\sigma_2\in M\setminus 2M$.

{\rm{(6)}} Assume $W_{\Pi}\cdot \Pi$ is $D_{l+1}^{(2)}$, and 
$k(\al_0)=2$, $g(\al_0)=2\bbbz +1$,  $k(\al_1)=1$,  $g(\al_1 )=2\bbbz +1$, $k(\beta)=1$,
$g(\beta)=\emptyset$  ($\beta\in\PIlg$),
so $\Lsh =M$, $\Llg =\{0,a\}+2M$, and $\Lex =a+2M$.
Then we have $\dim\fg_{\sigma_1}=l+1$ for all 
$\sigma_1\in \Llg\setminus\{0\}$, and $\dim\fg_{\sigma_2}=1$ for all 
$\sigma_2\in M\setminus \Llg$.
(This $R$ is isomorphic to $R(\Pi_2,k_2,g_2)$ for which
$W_{\Pi_2}\cdot \Pi_2$ is $A_{2l}^{(2)}$, 
and $k_2(\al )=1$, $g_2(\al )=\emptyset$
($\al\in\PIsh$), $k_2(\beta )=2$, $g_2(\beta )=\emptyset$
($\beta\in\PIlg\cup\PIex$).)

{\rm{(7)}} Assume $W_{\Pi}\cdot \Pi$ is $A_{2l}^{(2)}$, and $k(\al)=1$ ($\al\in\Pi$),
$g(\al_1)=2\bbbz +1$, 
$g(\beta)=\emptyset$  ($\beta\in\PIlg\cup \PIex$),
so $\Lsh =\Llg =M$, and $\Lex =\{\delta,\delta+a,a\}+2M$.
Then we have $\dim\fg_\sigma=l+1$ for all 
$\sigma\in M\setminus\{0\}$.

{\rm{(8)}} Assume $W_{\Pi}\cdot \Pi$ is $B_l^{(1)}$, and $k(\al)=1$ ($\al\in\Pi$),
$g(\al_1)=2\bbbz +1$, 
$g(\beta)=\emptyset$  ($\beta\in\PIlg$),
so $\Lsh =\Llg =M$, and $\Lex =a+2M$.
Let $M^\prime=\{0,a\}+2M$.
Then we have $\dim\fg_{\sigma_1}=l+1$ for all 
$\sigma_1\in M^\prime\setminus\{0\}$, and $\dim\fg_{\sigma_2}=1$ for all 
$\sigma_2\in M\setminus M^\prime$.
(This $R$ is isomorphic to $R(\Pi_3,k_3,g_3)$ for which
$W_{\Pi_3}\cdot \Pi_3$ is $A_{2l}^{(2)}$, 
and $k_3(\al )=1$, $g_3(\al )=\emptyset$
($\al\in\PIsh\cup\PIlg$), $k_3(\beta )=2$, $g_3(\beta )=\emptyset$
($\beta\in\PIex$).)

{\rm{(9)}} Assume $W_{\Pi}\cdot \Pi$ is $A_{2l}^{(2)}$, and $k(\al)=1$, $g(\al)=\emptyset$ ($\al\in\Pi$),
so $\Lsh =\Llg =M$, and $\Lex =\{a,\delta+a\}+2M$.
Then we have $\dim\fg_{\sigma_1}=l+1$ for all 
$\sigma_1\in M^\prime\setminus(\Lex\cup\{0\})$, and $\dim\fg_{\sigma_2}=1$ for all 
$\sigma_2\in \Lex$.
(This $R$ is isomorphic to $R(\Pi_4,k_4,g_4)$ for which
$W_{\Pi_4}\cdot \Pi_4$ is $A_{2l}^{(2)}$, 
and $k_4(\al_1 )=1$, $g_4(\al_1 )=2\bbbz +1$,
$k_4(\al_0 )=2$, $g_4(\al_1 )=\emptyset$,
$k_4(\beta )=1$, $g_4(\beta )=\emptyset$
($\beta\in\PIlg$).)

{\rm{(10)}} Assume $W_{\Pi}\cdot \Pi$ is $D_{l+1}^{(2)}$, and $k(\al)=1$ ($\al\in\Pi$),
$g(\al_0)=2\bbbz +1$, $g(\beta)=\emptyset$ ($\beta\in\PIlg\cup\{\al_1\}$).
Then 
\begin{equation}\label{eqn:tenL}
\Lsh=M,\,\,\Llg=\{0,a\}+2M\,\,\mbox{and}\,\,\Lex=\{2\delta+a,2\delta+3a\}+4M,
\end{equation}
and we have
\begin{equation}\label{eqn:mlexx}
\dim\fg_{p\delta+za}=\left\{\begin{array}{ll}
l+1 & \mbox{if $p\equiv_4 0$ and ${p \brack z}\ne {0 \brack 0}$,} \\
1 & \mbox{if $p\equiv_2 1$,} \\
l & \mbox{if $p\equiv_4 2$ and $z\equiv_2 0$,} \\
l+1 & \mbox{if $p\equiv_4 2$ and $z\equiv_2 1$.}
\end{array}\right.
\end{equation}
(This $R$ is isomorphic to $R(\Pi_5,k_5,g_5)$ for which
$W_{\Pi_5}\cdot \Pi_5$ is $A_{2l}^{(2)}$, 
$k(\al)=(\al,\al)/(\gamma_1,\gamma_1)$, $g_5(\al )=\emptyset$
($\al\in\PIlg$).)

\end{thm}
%
%
%\subsection{Special example}
%
%In this subsection, we assume
%that $l\geq 2$, $W_\Pi.\Pi$ is $D_{l+1}^{(2)}$, $k(\al_i)=1$ for $0\leq i\leq l$,
%and $g(\al_0)=2\bbbz+1$, $g(\al_1)=\emptyset$ and $g(\al_j)=\emptyset$
%for $2\leq j\leq l-1$.
%Then $\Lsh=M$,
%$\Llg=2\bbbz\delta\oplus\bbbz a=\{0,a\}+2M$, 
%and $\Lex=(4\bbbz+2)\delta\oplus(2\bbbz +1)a=\{2\delta +a,2\delta +3a\}+4M$. 
%Let $\Pi_1\cup\{a_1\}$ be a fundamental-set of $R$. Let $\delta_1=\delta(\Pi_1)$, so
%$\{\delta_1,a_1\}$ is a $\bbbz$-basis of $M$.

\begin{figure}
  \begin{center}
\setlength{\unitlength}{1mm}
\begin{picture}(110,100)(-10,-10)
\put(0,0){$l+2$} \multiput(10,0)(10,0){9}{$l+1$}
\multiput(0,10)(10,0){10}{$1$}
\multiput(0,20)(20,0){5}{$l$}\multiput(10,20)(20,0){5}{$l+1$}
\multiput(0,30)(10,0){10}{$1$}
\multiput(0,40)(10,0){10}{$l+1$}
\multiput(0,50)(10,0){10}{$1$}
\multiput(0,60)(20,0){5}{$l$}\multiput(10,60)(20,0){5}{$l+1$}
\multiput(0,70)(10,0){10}{$1$}
\multiput(0,80)(10,0){10}{$l+1$}

\put(0,0){\vector(1,1){90}}
\put(0,0){\vector(1,2){45}}

\put(0,-5){${0 \brack 0}$}\put(7,-9){\vector(1,0){10}} \put(10,-6){$a$}
\put(-7,0){\vector(0,1){10}}\put(-11,2){$\delta$}
\end{picture} 
  \end{center}
  \caption{$\dim\fg_{m\delta +ra}$ in (\ref{eqn:mlexx})}
  \label{fig:figmlexx}
\end{figure}

(At this moment, we do not see why $\dim\fg_{p\delta+za}$ are periodic modulo $tM$ for some $t\in\bbbn$. 
Maybe one of reasons is that $\fg$ may be realized as a `fixed point' Lie algebra, 
see also \cite{ABY}, \cite{yoshii}.) 
\newline\par
{\it {Proof.}} We only prove (10), since (1)-(9) are similarly treated.

Assume $(\al_1,\al_1)=1$. Define $\varepsilon_i\in \v$ ($1\leq i\leq l$) by
$\varepsilon_1:=\al_1$ and $\varepsilon_j:=\al_j+\varepsilon_{j-1}$ ($2\leq j\leq l$).
Then $(\varepsilon_i,\varepsilon_j)=\delta_{ij}$, and $\al_0=\delta-\varepsilon_1$.
Moreover, we have
\begin{equation}\begin{array}{l}
W_\Pi\cdot \al_1=\cup_{\epsilon\in\{-1,1\},1\leq i\leq l}\epsilon\varepsilon_i+2\bbbz \delta, \\
W_\Pi\cdot \al_r=\cup_{\epsilon_1,\epsilon_2\in\{-1,1\},1\leq i<j\leq l}\epsilon_1\varepsilon_i+
\epsilon_2\varepsilon_j+2\bbbz \delta\,\,(2\leq r\leq l), \\
W_\Pi\cdot \al_0=\cup_{\epsilon\in\{-1,1\},1\leq i\leq l}\epsilon\varepsilon_i+(2\bbbz +1)\delta.
\end{array}
\end{equation} Then by (\ref{eqn:eqnexppre}), we have
\begin{equation}\begin{array}{l}
R= \cup_{\epsilon\in\{-1,1\},1\leq i\leq l}\epsilon\varepsilon_i+2\bbbz \delta +\bbbz a \\
\quad\quad \cup\cup_{\epsilon_1,\epsilon_2\in\{-1,1\},1\leq i<j\leq l}\epsilon_1\varepsilon_i+
\epsilon_2\varepsilon_j+2\bbbz \delta +\bbbz 2a \\
\quad\quad \cup\cup_{\epsilon\in\{-1,1\},1\leq i\leq l}\epsilon\varepsilon_i+(2\bbbz +1) \delta +\bbbz a \\
\quad\quad \cup\cup_{\epsilon\in\{-1,1\},1\leq i\leq l}2(\epsilon\varepsilon_i+(2\bbbz +1) \delta) +(2\bbbz +1)a \\
\quad= \cup_{\epsilon\in\{-1,1\},1\leq i\leq l}\epsilon\varepsilon_i+M \\
\quad\quad \cup\cup_{\epsilon_1,\epsilon_2\in\{-1,1\},1\leq i<j\leq l}\epsilon_1\varepsilon_i+
\epsilon_2\varepsilon_j+2M \\
\quad\quad \cup\cup_{\epsilon\in\{-1,1\},1\leq i\leq l}\epsilon2\varepsilon_i+(4\bbbz +2) \delta +(2\bbbz +1)a.
\end{array}
\end{equation} Hence we have (\ref{eqn:tenL}), as desired.

Let $\Pi^\prime\cup\{a^\prime\}$ be a fundamental-set of $R$. Let $\delta^\prime:=\delta(\Pi^\prime)$, so
$\{\delta^\prime,a^\prime\}$ is a $\bbbz$-basis of $M$.

Assume $a^\prime\equiv_4 a={0 \brack 1}$.
Then $\delta^\prime\equiv_4 \delta={1 \brack y}$, where we replace $\Pi^\prime$ with 
$-\Pi^\prime$ if necessary.
Let $\delta^{\prime\prime}=\delta^\prime-ya^\prime$. Then 
$\{\delta^{\prime\prime},a^\prime\}$ is a $\bbbz$-basis of $M$.
Since $\delta^{\prime\prime}\equiv_4 \delta={1 \brack 0}\equiv_2 \delta={1 \brack 0}$,
we have $\Llg=\{0,a^\prime\}+2M$ and 
$\Lex=\{2\delta^{\prime\prime} +a^\prime,2\delta^{\prime\prime} +3a^\prime\}+4M$.
Hence we have the root system isomorphism $f_1:\bbbr R\to\bbbr R$ (cf. (\ref{eqn:rtiso}))
such that $f_1(\al_j)=\al_j$ ($1\leq j\leq l$), $f_1(\delta)=\delta^{\prime\prime}$
and $f_1(a)=a^\prime$. Then by Theorem~\ref{theorem:mainuni}, we have 
$\dim\fg_{ma^\prime}=l+1$ for $m\in\bbbzt$.

Assume $a^\prime\equiv_4 \delta={1 \brack 0}$.
Let $R_5=R(\Pi_5,k_5,g_5)$ be as in the statement.
Let $\fg^\prime:=\fg (\Pi_5,k_5,g_5)$.
%Let $R_1=R(\Pi_1,k_1,g_1)$ be such that $W_{\Pi_1}.\Pi_1$ is $A_{2l}^{(2)}$, and 
%that $k(\al_1)=1$, 
%$k(\al_i)=2$ ($2\leq i\leq l$), $k(\al_0)=4$, and $g(\al_j)=\emptyset$
%($0\leq j\leq l$). 
Define the $\bbbr$-linear isometry
$f_2:\bbbr R_5 \to \bbbr R$ by 
$f_2(\al_j)=\al_j$ ($1\leq j\leq l$), $f_2(\delta)=2\delta-a$
and $f_2(a)=\delta$. Note that 
$f_2(\Lsh)=f_2(M)=M=\Lsh$, $f_2(\Llg)=f_2(\{0,\delta\}+2M)=\Llg$
and $f_2(\Lex)=f_2(\{\delta,3\delta\}+4M)=\Llg$. Hence $f_2$ is 
a root system isomorphism. 
Let $a^{\prime\prime}:=f_2^{-1}(a^\prime)$. 
Then $a^{\prime\prime}\equiv_4 a$.
By the same argument
as above, as for $\dim \fg^\prime _{ma^{\prime\prime}}$,
we have the same equalities as in (\ref{eqn:mlexpre}) below.
%Assume $a_1\equiv_4 \delta={1 \brack 0}$. Let 
%$C=\left[\begin{array}{cc} 2 & 1 \\ -1 & 0 \end{array}\right]$,  so $\det C=1$
%and $C^{-1}=\left[\begin{array}{cc} 0 & -1 \\ 1 & 2 \end{array}\right]$.
%Then $\delta_1=C{x \brack y}$ for some $x$, $y\in\bbbz$.
%Since $a_1\equiv_4 C{0 \brack 1}$, we have $\epsilon x\equiv_4 1$
%for some $\epsilon\in\{-1,1\}$. 
%Let $b=\epsilon(\delta_1-y a_1)$, so $\{b,a_1\}$ is a $\bbbz$-basis of $M$.
%Then $b\equiv_4 C{1 \brack 0}$. 
%Since $[b,a_1]\equiv_4 [\delta,a]C$, we have $\delta\equiv_4 a_1$
%and $a\equiv_4 -b+2a_1\equiv_2 b$.
%Hence $\Llg=\{0,b\}+2M=\bbbz b\oplus 2\bbbz a_1$
%and
%$\Lex=\{3b,b\}+4M=(2\bbbz +1)b\oplus 4\bbbz a_1$.
%Hence we conclude that $W_{\Pi_1}.\Pi_1$ is $A_{2l}^{(2)}$, and 
%that $k(\al_1)=1$, 
%$k(\al_i)=2$ ($2\leq i\leq l$), $k(\al_0)=4$, and $g(\al_j)=\emptyset$
%($0\leq j\leq l$). 
Then Theorem~\ref{theorem:mainuni} implies that
\begin{equation}\label{eqn:mlexpre}
\dim\fg_{ma^\prime}=\left\{\begin{array}{ll}
l+1 & \mbox{if $m\ne 0$ and $m\equiv_4 0$,} \\
1 & \mbox{if $m\equiv_2 1$,} \\
l & \mbox{if $m\equiv_4 2$.}
\end{array}\right.
\end{equation}

For other $a^\prime$'s,
we can utilize %Theorem~\ref{theorem:mainuni} and
the %three 
root system isomorphisms $f_i:\bbbz R\to \bbbz R$
($3\leq i\leq 5$)  defined by
$f_i(\al_j)=\al_j$ for all $1\leq j\leq l$, and 
$f_3({1 \brack 0})={-1 \brack 0}$, $f_3({0 \brack 1})={0 \brack 1}$,  
$f_4({1 \brack 0})={1 \brack 0}$, $f_4({0 \brack 1})={0 \brack -1}$,   
$f_5({1 \brack 0})={1 \brack 1}$, $f_5({0 \brack 1})={0 \brack 1}$.
%$f_4({1 \brack 0})={1 \brack -1}$, $f_4({0 \brack 1})={0 \brack 1}$.  
Let $R_6=R(\Pi_6,k_6,g_6)$ be such that $W_{\Pi_6}\cdot \Pi_6$ is $D_{l+1}^{(2)}$,
$k_6(\al_i)=1$ for $1\leq i\leq l$,
and $g_6(\al_0)=\emptyset$, $g_6(\al_1)=2\bbbz+1$ and $g_6(\al_j)=\emptyset$
for $2\leq j\leq l-1$. Then we can also use 
the root system isomorphism $f_6:\bbbz R_6\to \bbbz R$
defined by
$f_6(\al_j)=\al_j$ ($1\leq j\leq l$), 
$f_6(\delta)=\delta$ and $f_6(a)=2\delta +a$.
%the two root system isomorphisms $f_i:\bbbz R_2\to \bbbz R$
%($5\leq i\leq 6$) defined by
%$f_i(\al_0)=\al_1$, $f_i(\al_1)=\al_0$,
%$f_i(\al_j)=\al_{l-j+1}$ for all $2\leq j\leq l-1$, 
%$f_i({1 \brack 0})={1 \brack 0}$, $f_5({0 \brack 1})={2 \brack 1}$, 
%$f_6({0 \brack 1})={-2 \brack 1}$.

Finally we have
\newline\par
Case-1. If $a^\prime\equiv_4 {0 \brack 1}$, ${0 \brack 3}$, ${2 \brack 1}$
or ${2 \brack 3}$, then we have $\dim\fg_{ma^\prime}=l+1$
for $m\in\bbbzt$.

Case-2. If $a^\prime\equiv_4 {1 \brack 0}$, ${3 \brack 0}$, 
${1 \brack 1}$, ${1 \brack 3}$, ${3 \brack 1}$, ${3 \brack 3}$, 
${1 \brack 2}$ or ${3 \brack 2}$, 
then the same as (\ref{eqn:mlexpre}) holds.
\newline\par Let $\lambda=p\delta+za={p \brack z}=m a^\prime$ with $p$, $z\in\bbbz$
and $m\in\bbbzt$.
Let ${x \brack y}=a^\prime$, so $x\bbbz+y\bbbz=\bbbz$.

Assume that $p\equiv_4 0$. If $x\equiv_2 1$, then $m\equiv_4 0$, so $\dim\fg_\lambda=l+1$.
If $x\equiv_2 0$, then $y\equiv_2 1$, so Case-1 implies $\dim\fg_\lambda=l+1$.

Assume that $p\equiv_4 2$ and $z\equiv_2 0$.
If $x \equiv_2 0$, then $y \equiv_2 1$, so $m \equiv_2 0$, so $p\equiv_4 0$, contradiction.
Hence $x\equiv_2 1$, so $m\equiv_4 2$, so Case-2 implies $\dim\fg_\lambda=l$.

Assume that $p\equiv_4 2$ and $z\equiv_2 1$.
Then $m\equiv_2 1$, $y\equiv_2 1$ and $x\equiv_2 0$, so Case-1 implies $\dim\fg_\lambda=l+1$.

Assume that $p\equiv_2 1$.
Then $m\equiv_2 1$ and $x\equiv_2 1$, so Case-2 implies $\dim\fg_\lambda=1$.

Thus we have (\ref{eqn:mlexx}), as desired. This completes the proof. 
\hfill $\Box$
%Thus we have 
%\begin{equation}\label{eqn:mlex}
%\dim\fg_{p\delta+za}=\left\{\begin{array}{ll}
%l+1 & \mbox{if $p\equiv_4 0$ and ${p \brack z}\ne {0 \brack 0}$,} \\
%1 & \mbox{if $p\equiv_2 1$,} \\
%l & \mbox{if $p\equiv_4 2$ and $z\equiv_2 0$,} \\
%l+1 & \mbox{if $p\equiv_4 2$ and $z\equiv_2 1$.}
%\end{array}\right.
%\end{equation}

\vspace{1cm}

{\bf{Acknowledgment.}} The authors would like to express their thanks to Kaiming Zhao for a fruitful  discussion on 
Remark~\ref{remark:mod}. Most of this paper is motivated by delivery of talks and lectures
by the authors at the University of Isfahan, Qom University and Arak 
University in Iran in September 2005 and December 2008. The authors would like to express their thanks to
professors, staffs and students there
for hospitality, and careful and warm listening to the authors' 
presentations.

\newpage
{\small \sc Saeid Azam,
Department of Mathematics, the University of Isfahan, Isfahan, Iran,
P.O.Box 81745-163} 

{\small \textit{E-mail address:} \texttt{azam@sci.ui.ac.ir}}

\vspace{.5\baselineskip}

{\small \sc Hiroyuki Yamane,
Department of Pure and Applied Mathematics,
Graduate School of Information Science
and Technology, Osaka University, Toyonaka 560-0043,
Japan}

{\small \textit{E-mail address:} \texttt{yamane@ist.osaka-u.ac.jp}}   

\vspace{.5\baselineskip}
{\small \sc Malihe Yousofzadeh, Department of Mathematics, the University of Isfahan, Isfahan, Iran,
P.O.Box 81745-163}

{\small \textit{E-mail address:} \texttt{ma.yousofzadeh@sci.ui.ac.ir}}


\begin{thebibliography}{99}
\bibitem[1]
{Allisonetal}
B.~Allison, S.~Azam, S.~Berman, Y.~Gao,
A.~Pianzola,
Extended affine Lie algebras and their root systems,
\textit{Mem. Amer. Math. Soc.} \textbf{603}  (1997)
1--122

\bibitem[2]
{Azam} S.~Azam,
Extended affine root systems,
\textit{J.~Lie Theory} \textbf{12} (2002) 515--527

\bibitem[3]{ABY}
S.~Azam, S.~Berman, N.~Yousofzadeh, 
Fixed point subalgebras of extended affine Lie algebras. 
\textit{J. Algebra} \textbf{287} (2005), no. 2, 351--380 

\bibitem[4]{AKY} S.~Azam, V.~Khalili and M.~Yousofzadeh,  Extended affine root system of type BC, 
J. Lie theory \textbf{15} (2005). no.~1, 145-181

\bibitem[5]{AYY} S.~Azam, H.~Yamane, M.~Yousofzadeh, 
Universal coverings of Lie tori (A finite presentation),
Preprint, arXiv:0906.0158 

\bibitem[6]
%BGK]
{BermanGaoKrylyuk} S.~Berman, Y.~Gao,
Y.~Krylyuk,
Quantum tori and the structure of elliptic quasi-simple
Lie algebras,
\textit{J.~Funct. Anal.} \textbf{135} (1996) 339--389


\bibitem[7]{Drinfeld} V.G.~Drinfeld,
A new realization of Yangians and quantized affine
algebras,
\textit{Soviet Math. Dokl.} \textbf{36} (1988) 212-216

\bibitem[8]{HecYam08b}
I. Heckenberger, H. Yamane, A generalization of Coxeter groups, root systems, and Matsumoto's theorem, 
\textit{Math. Z.}, \textbf{259} (2008), 255-276

\bibitem[9]{Humphreys} J.E.~Humphreys,
Introduction to Lie algebras and
representation theory,
Graduate Texts in Mathematics, Vol.~9,
Springer-Verlag, New York Berlin, 1972


\bibitem[10]
{kac1} V.G.~Kac, Infinite dimensional
Lie algebras
(3rd ed.),  Cambridge Univ. Press, Cambridge, 1990

\bibitem[11]
{kac2} \bysame,
Infinite-dimensional algebras, Dedekind's
$\eta$-function,
classical M\"{o}bius function and the very strange
formula,
\textit{Adv. in Math.} \textbf{30} (1978) 85--136


\bibitem[12]{lang} S.~Lang, {\it{Algebra}}, Addison-Weslay Publishing 
Company Inc. 1965 


\bibitem[13]
{MacDonald} I.G.~MacDonald, Affine root
systems
and Dedekind's $\eta$-function,
\textit{Invent. Math.} \textbf{15} (1972) 91--143


\bibitem[14]
{MacDonald2} \bysame, {\it{Affine Hecke algebras and orthogonal polynomials}},
Cambridge tracts in mathematics 157, 
Cambridge university press, 2003

\bibitem[15]
{MoodyRaoYokonuma} R.V.~Moody,
S.~Eswara~Rao,
T.~Yokonuma,
Toroidal Lie algebras and vertex representations,
\textit{Geom. Dedicata} \textbf{35} (1990) 283--307


\bibitem[16]
{saito1}
K.~Saito, Extended affine root systems. I. Coxeter
transformations,
\textit{Publ. Res. Math. Inst. Sci.} \textbf{21} (1985)
75--179

\bibitem[17]
{saitotakebayashi} K.~Saito, T.~Takebayashi,
Extended affine root systems.
III.
Elliptic
Weyl groups,
\textit{Publ. Res. Math. Inst. Sci.} \textbf{33} (1997)
301--329

\bibitem[18]
{saitoyoshii} K.~Saito, D.~Yoshii,
Extended affine root
system. IV.
Simply-laced elliptic Lie algebras,
\textit{Publ. Res. Math. Inst. Sci.} \textbf{36} (2000)
385--421


\bibitem[19]
{yamane2} H.~Yamane,
A Serre-type theorem for the elliptic Lie algebras with
rank $\geq 2$, \textit{Publ. Res. Math. Inst. Sci.}
\textbf{40} (2004) 441--469

\bibitem[20]{yoshii} D.~Yoshii, Elliptic Lie algebras
(Inhomogeneous
Cases),
{\it{Preprint}}, RIMS-1308 (2001)

\bibitem[21]{Yousofzadeh} M.~Yousofzadeh, 
A generalization of extended affine Lie algebras, Comm. Algebra 
\textbf{35} (2007), no 12, 4277-4302 

\bibitem[22]{YouRIMS} M. Yousofzadeh,
A presentation of Lie tori of type $B_\ell$, Publ. Res. Inst. Math. Sci.  {\bf 44}  (2008),  no. 1, 1--44





\end{thebibliography}
\end{document}